\documentclass[11pt]{article}
\usepackage{amssymb,amsfonts,amsthm}

\newcounter{ppp}
\newcommand{\be}[1]{\begin{equation}\label{#1}}
\newcommand{\ee}{\end{equation}}

\newcommand{\topp}[1]{\left\lceil{#1}\right\rceil}
\newcommand{\bott}[1]{\left\lfloor{#1}\right\rfloor}
\newcommand{\suff}{\hbox{\bf suff}}

\newcommand{\mm}{{\cal M}}

\newcommand{\nn}{{\cal N}}
\newcommand{\wtk}{\widetilde\kk}
\newcommand{\wtD}{\widetilde\Delta}
\newcommand{\wtF}{\widetilde\F^+}
\newcommand{\td}{{\tilde\delta}}

\newcommand{\tf}{{\tilde f}}
\newcommand{\tp}{{\tilde p}}
\newcommand{\tq}{{\tilde q}}
\newcommand{\tr}{{\tilde r}}
\newcommand{\ts}{{\tilde s}}
\newcommand{\tu}{{\tilde u}}
\newcommand{\tv}{{\tilde v}}

\newcommand{\dd}{{\cal D}}

\newcommand{\tos}{\stackrel{+}{\to}}
\newcommand{\kk}{{\cal K}}

\newcommand{\cd}{\mathop{\rm cd}}
\newcommand{\gd}{\mathop{\rm gd}}
\newcommand{\pp}{{\cal P}}

\newcommand{\uu}{{\cal U}}
\newcommand{\qq}{{\cal Q}}

\newcommand{\fft}{{\cal V}}
\newcommand{\hh}{{\cal H}}
\newcommand{\sss}{{\cal S}}
\newcommand{\fff}{{\cal F}}
\newcommand{\Ker}{\mathop{\rm Ker}}
\newcommand{\Img}{\mathop{\rm Im}}
\newcommand{\Ev}{\phi}
\newcommand{\Sqt}{\mathop{\rm\widetilde{Sq}}}
\newcommand{\zz}{{\mathbb{Z}}}
\newcommand{\ve}{\varepsilon}

\renewcommand{\wr}{\mathrel{\rm wr}}
\newcommand{\la}{\langle\,}
\newcommand{\ra}{\,\rangle}
\newcommand{\Sq}{\mathop{\rm Sq}}

\newcommand{\E}{{\bf E}}
\newcommand{\F}{{\bf F}}
\newcommand{\PP}{{\bf P}}

\begin{document}

\title
{Diagram groups and directed $2$-complexes: homotopy and homology}
\author{V.\,S.\,Guba,\ \,M.\,V.\,Sapir\thanks{The research of the first
author was partially supported by the RFFI grant 99--01--00894 and the
INTAS grant 99--1224. The research of the second author was
supported in part by the NSF grants DMS 0072307 and 9978802 and
the US-Israeli BSF grant 1999298.}}
\date{}

\maketitle

\begin{abstract}
We show that diagram groups can be viewed as fundamental groups of
spaces of positive paths on directed $2$-complexes (these spaces of
paths turn out to be classifying spaces). Thus diagram groups are
analogs of second homotopy groups, although diagram groups are as
a rule non-Abelian. Part of the paper is a review of the previous
results from this point of view. In particular, we show that the
so called rigidity of the R.\,Thompson's group $F$ and some other
groups is similar to the flat torus theorem. We find several
finitely presented diagram groups (even of type ${\cal F}_\infty$)
each of which contains all countable diagram groups. We show how
to compute minimal presentations and homology groups of a large
class of diagram groups.  We show that the Poincar\'e series of
these groups are rational functions. We prove that all integer
homology groups of all diagram groups are free Abelian.
\end{abstract}

\theoremstyle{plain}
\newtheorem{thm}{Theorem}[section]
\newtheorem{lm}[thm]{Lemma}
\newtheorem{cy}[thm]{Corollary}
\newtheorem{df}[thm]{Definition}
\newtheorem{rk}[thm]{Remark}
\newtheorem{ex}[thm]{Example}
\newtheorem{prop}[thm]{Proposition}
\newtheorem{prob}[thm]{Problem}

\tableofcontents
\section{Introduction}
\label{Intro}

The first definition of diagram groups was given by Meakin and
Sapir in terms of string rewriting systems (semigroup presentations).
Some results about diagram groups were obtained by Meakin's student
Vesna Kilibarda \cite{KilDiss,Kil}. Further results about diagram
groups have been obtained by the authors of this paper
\cite{GuSa97,GuSa99, GuSa02}, D.\,Farley \cite{Far00}, and
B.\,Wiest \cite{Wi02}.

The definition of diagram groups in terms of string rewriting
systems does not reflect the geometry of diagram groups and
geometrical nature of the constructions that can be applied to
diagram groups. In this paper we introduce a more geometric
definition of diagram groups in terms of directed $2$-complexes.

A directed $2$-complex (see \cite{Ste93,PW}) is a directed graph
equipped with $2$-cells each of which is bounded by two directed
paths (the top path and the bottom path). With any directed
$2$-complex one can associate the set of (directed) homotopies or
$2$-paths which is defined similar to the set of combinatorial
homotopies between $1$-paths in ordinary combinatorial $2$-complexes.
Equivalence classes of $2$-paths form a groupoid with respect to the
natural concatenation of homotopies. The local groups of that
groupoid are the diagram groups of the directed $2$-comp\-lex. Thus,
from this point of view, the diagram groups are ``directed"
analogs of the second homotopy groups.

The new point of view gave us an opportunity to revisit some
earlier results about diagram groups. We show that (a
multi-dimensional version of) the Squier complex of a semigroup
presentation has a natural realization as the space of positive
paths in a directed $2$-complex. We also show that several facts
about diagram groups proved earlier have a natural topological
interpretation in terms of directed $2$-complexes.

The only major part of our previous work that is not revisited
here is the combinatorics on diagrams (\cite[Section 15]{GuSa97}):
the description of centralizers of elements in diagram groups,
solving conjugacy problem, etc. All these results can be easily
generalized to the case of diagrams over directed $2$-complexes.

The largest part of the paper is devoted to completely new
results. In particular, we find several diagram groups of type
$\fff_\infty$ each of which contains all countable diagram groups.
One of them has only 3 generators and 6 defining relations. Recall
that a group $G$ is said to be of type $\fff_\infty$ if there is
some $K(G,1)$ CW complex having a finite $n$-skeleton in each
dimension $n$.

We show that the universal cover of the space of positive paths of a
directed $2$-complex $\kk$ is homeomorphic to the the space of positive
paths of what we call the {\em universal $2$-cover\/} of $\kk$ which is
again a directed $2$-complex. The universal $2$-covers possess some
remarkable properties which make them $2$-dimensional analogs of rooted
trees. We call such directed $2$-complexes {\em rooted $2$-trees\/}.

We study complete directed $2$-complexes (they are analogs of complete
string rewriting systems). We show how to construct a minimal $K(G,1)$
CW complex (with respect to the number of cells in each dimension) for
a diagram group $G$ of a complete directed $2$-complex. We compute integer
homology of such groups, and show that in the case when the groups are of
type $\fff_\infty$ (that happens very often), the Poincar\'e series are
rational. We answer Pride's question by showing that the integer homology
groups of arbitrary diagram groups are free Abelian. We also study the
cohomological dimension of diagram groups of complete directed $2$-complexes.
In particular, we show that the cohomological dimension of a group in that
class is $\ge n$ if and only if the group contains a copy of $\zz^n$ (for
any natural number $n$).

It was shown by Farley \cite{Far00} that diagram groups of finite semigroup
presentations act freely cellularly by isometries on CAT(0) cubical
complexes. One of the important results in the theory of CAT(0) groups is
the flat torus theorem that shows a rigid connection between a group acting
``nicely" on a CAT(0) space, and a geometric property of the space. The
algebraic property is ``to contain a copy of $\zz^n$", and the geometric
property is ``to contain a $\zz^n$-invariant copy of $\mathbb{R}^n$". We
show that similar rigid connection exists in our situation between, say,
the R.\,Thompson group $F$ and the universal cover of the space of positive
paths of the Dunce hat.

The results of this paper are used in our next paper \cite{GuSa02b} to
show that all diagram groups are totally orderable.

\section{Combinatorial definition}
\label{comb}

We start by giving a precise definition of directed $2$-complexes.
Our definition differs insignificantly from the original definition in
\cite{Ste93} and is close to the definition of \cite{PW}.

\begin{df}
\label{dirc}
{\rm
For every directed graph $\Gamma$ let $\PP$ be the set of all (directed)
paths in $\Gamma$, including the empty paths. A {\em directed $2$-complex\/}
is a directed graph $\Gamma$ equipped with a set $\F$ (called the {\em set
of $2$-cells\/}), and three maps $\topp{\cdot}\colon\F\to\PP$,
$\bott{\cdot}\colon\F\to\PP$, and $^{-1}\colon\F\to\F$ called {\em top\/},
{\em bottom\/}, and {\em inverse\/} such that
\begin{itemize}
\item for every $f\in\F$, the paths $\topp{f}$ and $\bott{f}$ are non-empty
and have common initial vertices and common terminal vertices,
\item $^{-1}$ is an involution without fixed points, and
$\topp{f^{-1}}=\bott{f}$, $\bott{f^{-1}}=\topp{f}$ for every $f\in\F$.
\end{itemize}
}
\end{df}

We shall often need an orientation on $\F$, that is, a subset
$\F^+\subseteq\F$ of {\em positive\/} $2$-cells, such that $\F$ is the
disjoint union of $\F^+$ and the set $\F^-=(\F^+)^{-1}$ (the latter is called
the set of {\em negative\/} $2$-cells).

If $\kk$ is a directed $2$-complex, then paths on $\kk$ are called
{\em $1$-paths\/} (we are going to have $2$-paths later). The initial and
terminal vertex of a $1$-path $p$ are denoted by $\iota(p)$ and
$\tau(p)$, respectively. For every $2$-cell $f\in\F$, the vertices
$\iota(\topp{f})=\iota(\bott{f})$ and $\tau(\topp{f})=\tau(\bott{f})$ are
denoted $\iota(f)$ and $\tau(f)$, respectively.

A $2$-cell $f\in\F$ with top $1$-path $p$ and bottom $1$-path $q$ will be
called a {\em $2$-cell of the form\/} $p=q$. We shall use a notation
$\kk=\la\E\mid\topp{f}=\bott{f},\ f\in\F^+\ra$ for a directed $2$-complex
with one vertex, the set of edges $\E$ and the set of $2$-cells $\F$.

For example, the complex $\la x\mid x^2=x\ra$ is the {\em Dunce hat\/}
obtained by identifying all edges in the triangle (Figure \theppp)
according to their directions.  It has one vertex, one edge, and one
positive $2$-cell. The remarkable feature of the Dunce hat is that the
famous R.\,Thompson's group $F$ is its diagram group (see Section
\ref{complete} below). (The survey \cite{CFP} collects some known results
about $F$. See also \cite{BrGe,BrSq85,Bro87,GuSa97,GuSa99,Brin,Gu01} for
other results about this group.)

\begin{center}
\unitlength=0.7mm
\special{em:linewidth 0.4pt}
\linethickness{0.4pt}
\begin{picture}(61.00,36.00)
\put(1.00,6.00){\vector(1,1){30.00}}
\put(31.00,36.00){\vector(1,-1){30.00}}
\put(1.00,6.00){\vector(1,0){60.00}}
\put(31.00,1.00){\makebox(0,0)[cc]{$x$}}
\put(14.00,24.00){\makebox(0,0)[cc]{$x$}}
\put(47.00,24.00){\makebox(0,0)[cc]{$x$}}
\end{picture}

\nopagebreak[4] Figure \theppp.
\end{center}
\addtocounter{ppp}{1}

There exists a natural way to assign a directed $2$-complex to every
semigroup presentation (to a string rewriting system). It is similar
to assigning a $2$-complex to any group presentation. If
$\pp=\la X\mid u_i=v_i\ \ (i\in I)\ra$ is a string rewriting system,
then the corresponding directed complex $\kk_\pp$ has one vertex, one
edge $e$ for each generator from $X$, and one positive $2$-cell for
each relation $u_i=v_i$. The top path of this cell is labelled by $u_i$
and the bottom path labelled by $v_i$.

Every directed $2$-complex $\kk=\la\E\mid\topp{f}=\bott{f},\ f\in\F^+\ra$
with one vertex can be considered as a rewriting system with the alphabet
$\E$ and the set of rewriting rules $\F^+$. The difference between these
directed $2$-complexes and string rewriting systems is that there may be
several $2$-cells in $\kk$ with the same top and bottom $1$-paths, hence
a rewriting rule $p=q$ can repeat many times. We shall observe later that
directed $2$-complexes with one vertex provide the same class of diagram
groups as all directed $2$-complexes. But sometimes it is convenient to
consider complexes with more than one vertex.

With the directed 2-complex $\kk$, one can associate a category $\Pi(\kk)$
whose objects are $1$-paths, and morphisms are {\em $2$-paths\/} (or
{\em homotopies\/}), i.\,e. sequences of replacements of $\topp{f}$ by
$\bott{f}$ in $1$-paths, $f\in\F$. More precisely, an {\em atomic $2$-path\/}
(an {\em elementary homotopy\/}) is a triple $(p,f,q)$, where $p$, $q$ are
$1$-paths in $\kk$, and $f\in\F$ such that $\tau(p)=\iota(f)$,
$\tau(f)=\iota(q)$. If $\delta$ is the atomic $2$-path $(p,f,q)$, then
$p\topp{f}q$ is denoted by $\topp{\delta}$, and $p\bott{f}q$ is denoted by
$\bott{\delta}$; these are called the {\em top\/} and the {\em bottom\/}
$1$-paths of the atomic $2$-path. Every {\em nontrivial\/} $2$-path $\delta$
on $\kk$ is a sequence of atomic paths $\delta_1$, \dots, $\delta_n$, where
$\bott{\delta_i}=\topp{\delta_{i+1}}$ for every $1\le i<n$. In this case $n$
is called the {\em length\/} of the $2$-path $\delta$. The {\em top\/} and
the {\em bottom\/} $1$-paths of $\delta$, denoted by $\topp{\delta}$ and
$\bott{\delta}$, are $\topp{\delta_1}$ and $\bott{\delta_n}$, respectively.
We say that $\delta$ is {\em positive\/} if each $\delta_i$ corresponds to a
positive $2$-cell in $\kk$. Every $1$-path $p$ is considered as a trivial
$2$-path with $\topp{p}=\bott{p}=p$. These are the identity morphisms in the
category $\Pi(\kk)$. The composition of $2$-paths $\delta$ and $\delta'$ is
called {\em concatenation\/} and is denoted $\delta\circ\delta'$.
\vspace{0.5ex}

We say that $1$-paths $p$, $q$ in $\kk$ are (directly) {\em homotopic\/}
whenever there exists a $2$-path $\delta$ such that $\topp{\delta}=p$ and
$\bott{\delta}=q$. We also say that $\delta$ {\em connects\/} $p$ to $q$ in
$\kk$.
\vspace{0.5ex}

As in the standard homotopy theory (see, for example, \cite{Si93}), we need to
introduce a homotopy relation on the set of $2$-paths and identify homotopic
$2$-paths in $\Pi(\kk)$. To do this, we choose a computation-friendly way,
similar to the one developed by Peiffer, Reidemeister and Whitehead for the
second homotopy group of a combinatorial $2$-complex, and later simplified by
Huebschmann, Sieradski and Fenn (see Bogley and Pride \cite{BP93}). The idea is
to represent the elements of the second homotopy groups in terms of the so
called pictures. We are going to use the dual objects called diagrams (for a
picture version, see \cite{GuSa97}).

With every atomic $2$-path $\delta=(p,f,q)$, where $\topp{f}=u$, $\bott{f}=v$
we associate the labelled plane graph $\Delta$ on Figure \theppp. An arc
labelled by a word $w$ is subdivided into $|w|$ edges\footnote{In this paper,
we denote the length of a word or a path $w$ by $|w|$.}. All edges are
oriented from the left to the right. The label of each oriented edge of the
graph is a symbol from the alphabet $\E$, the set of edges in $\kk$. As a
plane graph, it has only one bounded face; we label it by the corresponding
cell $f$ of $\kk$. This plane graph $\Delta$ is called the diagram of $\delta$.
Such diagrams are called {\em atomic\/}. The leftmost and rightmost vertices
of $\Delta$ are denoted by $\iota(\Delta)$ and $\tau(\Delta)$, respectively.
The diagram $\Delta$ has two distinguished paths from $\iota(\Delta)$ to
$\tau(\Delta)$ that correspond to the top and bottom paths of $\delta$. Their
labels are $puq$ and $pvq$, respectively. These are called the top and the
bottom paths of $\Delta$ denoted by $\topp{\Delta}$ and $\bott{\Delta}$.

\begin{center}
\unitlength=1.00mm
\special{em:linewidth 0.4pt}
\linethickness{0.4pt}
\begin{picture}(110.66,26.33)
\put(0.00,13.33){\circle*{1.33}}
\put(40.00,13.33){\circle*{1.33}}
\put(70.00,13.33){\circle*{1.33}}
\put(110.00,13.33){\circle*{1.33}}
\bezier{200}(40.00,13.33)(55.00,33.33)(70.00,13.33)
\bezier{200}(40.00,13.33)(55.00,-6.67)(70.00,13.33)
\put(20.00,17.33){\makebox(0,0)[cc]{$p$}}
\put(55.00,25.33){\makebox(0,0)[cc]{$u$}}
\put(90.00,17.33){\makebox(0,0)[cc]{$q$}}
\put(55.00,12.33){\makebox(0,0)[cc]{$f$}}
\put(55.00,0.00){\makebox(0,0)[cc]{$v$}}
\put(0.33,13.33){\line(1,0){39.67}}
\put(70.00,13.33){\line(1,0){40.00}}
\end{picture}
\vspace{1ex}

\nopagebreak[4] Figure \theppp.
\end{center}
\addtocounter{ppp}{1}

The diagram corresponding to the trivial $2$-path $p$ is just an arc labelled
by $p$; it is called a {\em trivial diagram\/} and it is denoted by $\ve(p)$.

Let $\delta=\delta_1\circ\delta_2\circ\cdots\circ\delta_n$ be a $2$-path
in $\kk$, where $\delta_1$, \dots, $\delta_n$ are atomic $2$-paths. Let
$\Delta_i$ be the atomic diagram corresponding to $\delta_i$. Then the
bottom path of $\Delta_i$ has the same label as the top path of $\Delta_{i+1}$
($1\le i<n$). Hence we can identify the bottom path of $\Delta_i$ with the
top path of $\Delta_{i+1}$ for all $1\le i<n$, and obtain a plane graph
$\Delta$, which is called the {\em diagram of the $2$-path\/} $\delta$.

It is clear that the above diagram $\Delta$, as a plane graph, has
exactly $n$ bounded faces or {\em cells\/}.

Two diagrams are considered {\em equal\/} if they are isotopic as plane
graphs. The isotopy must take vertices to vertices, edges to edges, it must
also preserve labels of edges and inner labels of cells. Two $2$-paths are
called {\em isotopic\/} if the corresponding diagrams are equal.

For example, consider the diagram on Figure \theppp\ below. It is clear
that it corresponds to the $2$-path
\be{sq1}
(p,f_1,qu_2r)\circ(pv_1q, f_2,r)
\ee
as well as the $2$-path
\be{sq2}
(pu_1q, f_2, r)\circ(p, f_1, qv_2r).
\ee

\begin{center}
\unitlength=1.00mm
\special{em:linewidth 0.4pt}
\linethickness{0.4pt}
\begin{picture}(122.66,24.00)
\put(2.00,13.00){\circle*{1.33}}
\put(27.00,13.00){\circle*{1.33}}
\put(47.00,13.00){\circle*{1.33}}
\put(77.00,13.00){\circle*{1.33}}
\put(97.00,13.00){\circle*{1.33}}
\put(122.00,13.00){\circle*{1.33}}
\bezier{152}(27.00,13.00)(37.00,30.00)(47.00,13.00)
\bezier{136}(27.00,13.00)(37.00,-1.00)(47.00,13.00)
\bezier{156}(77.00,13.00)(87.00,30.00)(97.00,13.00)
\bezier{144}(77.00,13.00)(87.00,-2.00)(97.00,13.00)
\put(15.00,15.00){\makebox(0,0)[cc]{$p$}}
\put(37.00,24.00){\makebox(0,0)[cc]{$u_1$}}
\put(37.00,14.00){\makebox(0,0)[cc]{$f_1$}}
\put(62.00,15.00){\makebox(0,0)[cc]{$q$}}
\put(87.00,24.00){\makebox(0,0)[cc]{$u_2$}}
\put(110.00,15.00){\makebox(0,0)[cc]{$r$}}
\put(37.00,3.00){\makebox(0,0)[cc]{$v_1$}}
\put(87.00,14.00){\makebox(0,0)[cc]{$f_2$}}
\put(87.00,3.00){\makebox(0,0)[cc]{$v_2$}}
\put(2.00,13.00){\line(1,0){25.00}}
\put(47.00,13.00){\line(1,0){30.00}}
\put(97.00,13.00){\line(1,0){25.00}}
\label{fig1}
\end{picture}

\nopagebreak[4] Figure \theppp.
\addtocounter{ppp}{1}
\end{center}

In that case we call the atomic $2$-paths $(p,f_1,qu_2r)$ and $(pu_1q,f_2,r)$
{\em independent\/}.

The proof of the following lemma in the case of semigroup presentations
can be extracted from the proof of \cite[Lemma 6.2]{GuSa97}. We leave it
to the reader to generalize the proof to the case of directed $2$-complexes.

\begin{lm}
\label{apis}
The isotopy relation on $2$-paths of a directed $2$-complex $\kk$ is the
smallest equivalence relation containing all pairs of $2$-paths of the
form {\rm (\ref{sq1})} and {\rm (\ref{sq2})} and invariant under concatenation.
\end{lm}

Concatenation of $2$-paths corresponds to the concatenation of diagrams:
if the bottom path of $\Delta_1$ and the top path of $\Delta_2$ have the
same labels, we can identify them and obtain a new diagram
$\Delta_1\circ\Delta_2$.

Note that for any atomic $2$-path $\delta=(p,f,q)$ in $\kk$ one can naturally
define its {\em inverse\/} $2$-path $\delta^{-1}=(p,f^{-1},q)$. The inverses of
all $2$-paths and diagrams are defined naturally. The inverse diagram $\Delta^{-1}$
of $\Delta$ is obtained by taking the mirror image of $\Delta$ with respect to a
horizontal line, and replacing labels of cells by their inverses.

Let us identify in the category $\Pi(\kk)$ all isotopic $2$-paths and also identify
each $2$-path of the form $\delta'\delta\delta^{-1}\delta''$ with $\delta'\delta''$.
The quotient category is obviously a groupoid (i.\,e. a category with invertible
morphisms). It is denoted by $\dd(\kk)$ and is called the {\em diagram groupoid\/}
of $\kk$. Two $2$-paths are called {\em homotopic\/} if they correspond to the same
morphism in $\dd(\kk)$. For each $1$-path $p$ of $\kk$, the local group of $\dd(\kk)$
at $p$ (i.\,e. the group of equivalence classes of $2$-paths connecting $p$ with
itself) is called the {\em diagram group of the directed $2$-complex $\kk$ with base\/}
$p$ and is denoted by $\dd(\kk,p)$. Notice that if $p$ is empty, then $\dd(\kk,p)$
is trivial by definition. In this paper, we shall usually ignore these diagram groups.

We shall give a much easier (equivalent) definition of the diagram groupoid
after Theorem \ref{thnf} below.

\begin{rk}
\label{rk33}
{\rm
Notice first that the diagram groups of a directed $2$-complex do not depend
on the orientation on the set of $2$-cells of that complex. Notice also that
if a directed $2$-complex $\kk'$ is obtained from $\kk$ by identifying vertices,
then the diagram groupoid of $\kk'$ may differ from the diagram groupoid of $\kk$
because the set of $1$-paths may increase, but the diagram groups of $\kk$ will
be diagram groups of $\kk'$ as well.
}
\end{rk}

One can easily check that if $\kk=\kk_\pp$ for some semigroup presentation $\pp$
and $w$ is a word over $X$ (that is, the corresponding path in $\kk_\pp$), then
the diagram group $\dd(\kk,w)$ we just defined coincides with the diagram group
$D(\pp,w)$ over $\pp$ defined in \cite{GuSa97}. Clearly, if $\kk=\kk_\pp$ then
$2$-paths are just the derivations over the semigroup presentation $\pp$.

It is convenient to define diagrams over a directed $2$-complex $\kk$
in an ``abstract" way, without referring to $2$-paths of $\kk$. Such
a definition was given by Kashintsev \cite{Kash70} and Remmers \cite{Rem}
in the case of semigroup presentations. Here we basically repeat their
definition and result.

\begin{df}
\label{diag}
{\rm
A diagram over $\kk=\la\E\mid\topp{f}=\bott{f},\,f\in\F^+\ra$ is
a finite plane directed and connected graph $\Delta$, where every
edge is labelled by an element from $\E$, and every bounded face
is labelled by an element of $\F$ such that:
\begin{itemize}
\item $\Delta$ has exactly one vertex-source $\iota$ (which has
no incoming edges) and exactly one vertex-sink $\tau$ (which has
no outgoing edges);
\item every $1$-path in $\Delta$ is simple;
\item each face of $\Delta$ labelled by $f\in\F$ is bounded by
two $1$-paths $u$ and $v$ such that the label of $u$ is $\topp{f}$,
the label of $v$ is $\bott{f}$, and the loop $uv^{-1}$ on the plane
goes around the face in the clockwise direction.
\end{itemize}
}
\end{df}

It is easy to see \cite{GuSa97} that every plane graph satisfying
the conditions of Definition \ref{diag} is situated between two
positive paths connecting $\iota$ and $\tau$. These paths are
$\topp{\Delta}$ and $\bott{\Delta}$.

We say that a diagram $\Delta$ over a directed $2$-complex $\kk$ is
a $(u,v)$-{\em diagram\/} whenever $u$ is the top label and $v$ is
the bottom label of $\Delta$. If $u$ and $v$ are the same, then
the diagram is called {\em spherical\/} (with base $u=v$).

The following lemma (see \cite[Lemma 3.5]{GuSa97}) shows that
diagrams over $\kk$ in the sense of Definition \ref{diag} are
exactly diagrams that correspond to $2$-paths in $\kk$. We will
often use this fact without reference.

\begin{lm}
\label{kash}
Let $\kk$ be a directed $2$-complex. Then $1$-paths $u$, $v$ are
homotopic in $\kk$ if and only if there exists a $(u,v)$-diagram over
$\kk$\ $($in the sense of Definition {\rm\ref{diag}}$)$.
\end{lm}

Diagrams over $\kk$ corresponding to homotopic $2$-paths are called
{\em equivalent\/}. The equivalence relation on the set of diagrams
(and the homotopy relation on the set of $2$-paths of $\kk$) can be
defined very easily as follows. We say that two cells $\pi_1$ and
$\pi_2$ in a diagram $\Delta$ over $\kk$ form a {\em dipole\/} if
$\bott{\pi_1}$ coincides with $\topp{\pi_2}$ and the labels of the
cells $\pi_1$ and $\pi_2$ are mutually inverse. Clearly, if $\pi_1$
and $\pi_2$ form a dipole, then one can remove the two cells from the
diagram and identify $\topp{\pi_1}$ with $\bott{\pi_2}$. The result
will be some diagram $\Delta'$. As in \cite{GuSa97}, it is easy to
prove that if $\delta$ is a $2$-path corresponding to $\Delta$, then
the diagram $\Delta'$ corresponds to a $2$-path $\delta'$, which is
homotopic to $\delta$. We call a diagram {\em reduced\/} if it does
not contain dipoles. A $2$-path $\delta$ in $\kk$ is called {\em reduced\/}
if the corresponding diagram is reduced. The following is an analog of
Kilibarda's lemma. The proof coincides with the proof of
\cite[Theorem 3.17]{GuSa97} and we omit it here.

\begin{thm}
\label{thnf}
Every equivalence class of diagrams over a directed $2$-complex
$\kk$ contains exactly one reduced diagram. Every $2$-path in $\kk$
is homotopic to a reduced $2$-path, every two homotopic reduced
$2$-paths have equal diagrams and so they contain the same number
of atomic factors.
\end{thm}

Thus one can define morphisms in the diagram groupoid $\dd(\kk)$ as reduced
diagrams over $\kk$ with operation ``con\-cat\-enat\-ion + reduc\-tion"
(that is, the product of two reduced diagrams $\Delta$ and $\Delta'$ is the
result of removing all dipoles from $\Delta\circ\Delta'$ step by step).

The diagram groupoid $\dd(\kk)$ has another natural operation, {\em addition\/}:
if $\Delta'$ and $\Delta''$ are diagrams over $\kk$ and
$\tau(\topp{\Delta'})=\iota(\topp{\Delta''})$ in $\kk$ then one can identify
$\tau(\Delta')$ with $\iota(\Delta'')$ to obtain the new diagram denoted by
$\Delta'+\Delta''$ and called the {\em sum\/} of $\Delta'$ and $\Delta''$. If
$\kk$ has only one vertex, then this operation is everywhere defined. In that
case the operation of addition makes $\dd(\kk)$ a tensor groupoid in the sense
of \cite{JS}.

Figure \theppp\ below illustrates the concepts of the concatenation of diagrams
and the sum of them.

\begin{center}
\unitlength=0.70mm
\special{em:linewidth 0.4pt}
\linethickness{0.4pt}
\begin{picture}(123.67,35.33)
\put(1.00,23.33){\circle*{1.33}}
\put(45.67,23.33){\circle*{1.33}}
\put(0.67,23.33){\line(1,0){45.00}}
\bezier{320}(0.67,23.33)(25.67,56.33)(45.67,23.33)
\bezier{332}(0.67,23.33)(24.67,-11.67)(45.67,23.33)
\put(23.67,30.33){\makebox(0,0)[cc]{$\Delta'$}}
\put(23.67,14.33){\makebox(0,0)[cc]{$\Delta''$}}
\put(23.67,0.67){\makebox(0,0)[cc]{$\Delta'\circ\Delta''$}}
\put(65.67,23.33){\circle*{1.33}}
\put(93.67,23.33){\circle*{1.33}}
\put(122.67,23.33){\circle*{1.33}}
\bezier{164}(65.67,23.33)(79.67,38.33)(93.67,23.33)
\bezier{152}(65.67,23.33)(80.67,10.33)(93.67,23.33)
\bezier{172}(93.67,23.33)(108.67,39.33)(122.67,23.33)
\bezier{168}(93.67,23.33)(109.67,8.33)(122.67,23.33)
\put(79.67,23.33){\makebox(0,0)[cc]{$\Delta'$}}
\put(108.67,23.33){\makebox(0,0)[cc]{$\Delta''$}}
\put(93.67,12.00){\makebox(0,0)[cc]{$\Delta'+\Delta''$}}
\end{picture}

\nopagebreak[4] Figure \theppp.
\end{center}
\addtocounter{ppp}{1}

\section{Topological definition}
\label{topdef}

We have seen that diagram groups are directed analogs of the
second homotopy groups. Recall that one can define the second
homotopy groups of a topological space as the fundamental group of
a space of paths. On the other hand, in the case of semigroup
presentations, the diagram groups can be defined as fundamental
groups of the so called $2$-dimensional Squier complexes associated
with the presentations (see \cite{GuSa97}).

In this section, we show that the Squier complex of a directed
$2$-complex $\kk$ can be considered as the space of certain positive
paths in $\kk$.

First we define a multi-dimensional version of the Squier complex $\Sq(\kk)$
from \cite{GuSa97} as a semi-cubical complex. Recall \cite{Ser51,BH} that a
{\em semi-cubical complex\/} is a family of disjoint sets $\{M_n;n\ge0\}$
(the elements of $M_n$ are called $n$-cubes) with face maps
$\lambda_{i}^k\colon M_n\to M_{n-1}$ ($1\le i\le n$, $k=0,1$) satisfying the
semi-cubical relations (the rightmost map act first):

\be{part}
\lambda^k_{i}\lambda^{k'}_{j}=\lambda^{k'}_{j-1}\lambda^k_{i}
\quad (i<j).
\ee

A {\em realization\/} of a semi-cubical complex $\{M_n;n\ge0\}$ can be
obtained as a factor-space of the disjoint union of Euclidean cubes, one
$n$-cube $c(x)$ for each element $x\in M_n$, $n\ge0$. The equivalence
relation identifies (point-wise) the cube $\lambda_{i}^k(c(x))$ with the
cube $c(\lambda_{i}^k(x))$ for all $i$, $k$, $x$. Here $\lambda_{i}^k(I^{n})$
is the corresponding $(n-1)$-face of the Euclidean $n$-cube $I^n$ (that is,
$\lambda_{i}^k(I^n)=I^{i-1}\times\{k\}\times I^{n-i}$ ).

\begin{df}
{\rm
The semi-cubical complex $\Sq(\kk)$ is defined as follows. For every
$n\ge0$, let $M_n$ be the set of {\em thin\/} diagrams \cite{Far00} of
the form
\be{thin}
\ve(u_0)+f^{(1)}+\ve(u_1)+\cdots+f^{(n)}+\ve(u_n)
\ee
where $f^{(i)}$ are negative\footnote{Taking negative edges instead of
positive simplifies some computations later.} $2$-cells of $\kk$
and $u_i$ are $1$-paths in $\kk$ (Figure 3 thus shows a thin diagram
with two cells). The face map $\lambda_{i}^k$ takes the thin diagram
$c$ of the form (\ref{thin}) to
\be{fop1}
\topp{c}_i=\ve(u_0)+f^{(1)}+\cdots+\topp{f^{(i)}}+\cdots+f^{(n)}+\ve(u_n)
\ee
if $k=0$ and
\be{fop2}
\bott{c}_i=\ve(u_0)+f^{(1)}+\cdots+\bott{f^{(i)}}+\cdots+f^{(n)}+\ve(u_n)
\ee
if $k=1$. These faces are called the {\em top\/} and the {\em bottom\/}
$i$th faces of $c$, respectively. It is easy to check that the conditions
(\ref{part}) are satisfied, so $\Sq(\kk)$ is a semi-cubical complex. }
\end{df}

\begin{rk}
\label{invar}
{\rm
Note that the realization of $\Sq(\kk)$ does not depend on the orientation
$\F^+$ on $\kk$. If we change the orientation, then some cells $f_i$ in the
thin diagram (\ref{thin}) are replaced by their inverses. This means that
we simply change the reference point of a cube. So we will think about every
thin diagram as of a cube with fixed reference point.
}
\end{rk}

Thus the vertices of $\Sq(\kk)$ are $1$-paths of $\kk$, the edges correspond
to negative atomic $2$-paths $(u,f,v)$, $2$-cells correspond to pairs of
independent atomic $2$-paths, etc. For example, the thin diagram
$\ve(u)+f+\ve(v)+g+\ve(w)$ determines a square with contour
$$
(u,f,v\topp{g}w)\circ(u\bott{f}v,g,w)\circ(u,f^{-1},v\bott{g}w)
\circ(u\topp{f}v,g^{-1},w).
$$

It is convenient to enrich the structure of the Squier complex
$\Sq(\kk)$ by introducing inverse edges: $(u,f,v)^{-1}=(u,f^{-1},v)$.
Then the edges $(u,f,v)$ will be called {\em positive\/} if $f$
is a positive $2$-cell of $\kk$, and {\em negative\/} if $f$ is
negative. As a result, the $1$-skeleton of $\Sq(\kk)$ turns into a
graph in the sense of Serre \cite{Se80}, and the $2$-skeleton of
$\Sq(\kk)$ coincides with the Squier complex defined in \cite{GuSa97}
provided $\kk=\kk_\pp$ for some $\pp$. Hence, in particular, the
fundamental groups of $\Sq(\kk_\pp)$ coincide with fundamental
groups of the Squier complex in \cite{GuSa97}.

Clearly the complex $\Sq(\kk)$ is in general disconnected. If $p$ is a
$1$-path in $\kk$, then by $\Sq(\kk,p)$ we will denote the connected
component of the Squier complex that contains $p$.

\begin{ex}
{\rm
Figure \theppp\ shows a part of the Squier complex of the Dunce hat on
Figure~1. The thick line shows the boundary of one of the $2$-cells in
this complex.

\begin{center}
\unitlength=1.00mm
\special{em:linewidth 0.4pt}
\linethickness{0.4pt}
\begin{picture}(67.00,24.33)
\put(7.00,13.00){\line(1,0){13.00}}
\put(20.00,13.00){\circle*{1.49}}
\put(7.00,13.00){\circle*{1.49}}
\put(33.00,13.00){\circle*{1.49}}
\put(46.00,13.00){\circle*{1.49}}
\put(59.00,13.00){\circle*{1.49}}
\bezier{88}(20.00,13.00)(27.00,21.67)(33.00,13.00)
\bezier{84}(20.00,13.00)(27.33,4.67)(33.33,13.00)
\linethickness{2.0pt}
\bezier{96}(33.00,13.00)(40.00,23.33)(46.00,13.00)
\bezier{88}(33.00,13.00)(40.33,4.00)(45.00,12.00)
\bezier{156}(46.50,15.00)(52.67,31.33)(59.00,13.00)
\bezier{160}(46.00,13.00)(53.33,-6.00)(59.00,13.00)
\linethickness{0.4pt}
\put(33.67,13.00){\line(1,0){12.67}}
\bezier{84}(46.33,13.00)(53.33,21.67)(59.00,13.00)
\bezier{88}(46.33,13.00)(54.00,4.00)(59.00,13.00)
\put(67.00,13.00){\makebox(0,0)[cc]{$\ldots$}}
\put(7.00,8.67){\makebox(0,0)[cc]{$x$}}
\put(20.00,8.67){\makebox(0,0)[cc]{$x^2$}}
\put(33.33,8.67){\makebox(0,0)[cc]{$x^3$}}
\put(45.66,8.67){\makebox(0,0)[cc]{$x^4$}}
\end{picture}

\nopagebreak[4] Figure \theppp.
\end{center}
\addtocounter{ppp}{1}
}
\end{ex}

The following theorem is similar to Kilibarda's statement (see
\cite[Theorem 6.1]{GuSa97}).

\begin{thm}
\label{fg}
Let $\kk$ be a directed $2$-complex, $p$ be a $1$-path in $\kk$.
Then the diagram group $\dd(\kk,p)$ is isomorphic to the
fundamental group $\pi_1(\Sq(\kk),p)$ of the semi-cubical complex
$\Sq(\kk)$ with the basepoint $p$.
\end{thm}

The proof of Kilibarda's theorem carries without any essential
changes. As an immediate corollary of Theorem \ref{fg}, we obtain
the following

\begin{cy}
\label{cy1}
Let $\kk$ be a directed $2$-complex, $p$ and $q$ be homotopic $1$-paths
in $\kk$. Then $\dd(\kk,p)$ is isomorphic to $\dd(\kk,q)$.
\end{cy}

\proof
Indeed, $p$ and $q$ belong to the same connected component of $\Sq(\kk)$.
\endproof

The diagram groups with different bases can be very
different but there exists the following  useful relationship
between them.

\begin{cy}
\label{cy2}
If $p=p_1p_2$ is a $1$-path in $\kk$, then $\dd(\kk,p_1)\times\dd(\kk,p_2)$
is embedded into $\dd(\kk,p_1p_2)$.
\end{cy}

\proof
Indeed, the map $(\Delta_1,\Delta_2)\to\Delta_1+\Delta_2$ from
$\dd(\kk,p_1)\times\dd(\kk,p_2)$ to $\dd(\kk,p_1p_2)$ is an
injective homomorphism (see \cite{GuSa97}, Remark 2 after Lemma 8.1).
\endproof

Now let us introduce a natural topological realization of the semi-cubical
complex $\Sq(\kk)$.

We expand the set of paths in $\kk$ allowing paths that go ``inside"
$2$-cells.

Let $\kk$ be a directed $2$-complex. Attaching a $2$-cell $f\in\F^+$ with
$p=\topp{f}$, $q=\bott{f}$ can be done as follows. Let $D=[0,1]\times[0,1]$
be a unit square. For any $t\in[0,1]$ we have the path $d_t$ in $D$ defined
by $d_t(s)=(s,t)\in D$ ($s\in[0,1]$). We attach this square to $\kk$ in such
a way that $d_0$ is identified with $p$, $d_1$ is identified with $q$, all
points of the form $(0,t)\in D$ are collapsed to $\iota(p)=\iota(q)$, all
points of the form $(1,t)$ are collapsed to $\tau(p)=\tau(q)$ ($t\in[0,1]$).

Now for any $t\in[0,1]$, the image of $d_t$ in $\kk$ becomes a path inside
the $2$-cell $f$. This path will be denoted by $f_t$. Clearly, $f_0=p$,
$f_1=q$. So we have a continuous family of paths $\{f_t\}$ ($t\in[0,1]$)
that transforms $p$ into $q$ (see Figure \theppp).

For any $f\in\F^+$, $t\in[0,1]$ one can also define $(f^{-1})_t=f_{1-t}$.
So $f_t$ makes sense for any $f\in\F$.

\begin{center}
\unitlength=1.00mm
\special{em:linewidth 0.4pt}
\linethickness{0.4pt}
\begin{picture}(113.66,52.00)
\put(3.00,36.00){\vector(1,0){40.00}}
\put(23.00,33.00){\makebox(0,0)[cc]{$d_t$}}
\put(-2.00,10.00){\makebox(0,0)[cc]{$(0,1)$}}
\put(49.00,10.00){\makebox(0,0)[cc]{$(1,1)$}}
\put(-2.00,50.00){\makebox(0,0)[cc]{$(0,0)$}}
\put(49.00,50.00){\makebox(0,0)[cc]{$(1,0)$}}
\put(49.00,36.00){\makebox(0,0)[cc]{$(1,t)$}}
\put(-2.00,36.00){\makebox(0,0)[cc]{$(0,t)$}}
\put(22.00,2.00){\makebox(0,0)[cc]{{\Large $D$}}}
\put(63.00,30.00){\circle*{1.33}}
\put(113.00,30.00){\circle*{1.33}}
\bezier{324}(63.00,30.00)(88.00,62.00)(113.00,30.00)
\bezier{308}(63.00,30.00)(88.00,1.00)(113.00,30.00)
\put(88.00,12.00){\makebox(0,0)[cc]{$q=f_1$}}
\put(88.00,49.00){\makebox(0,0)[cc]{$p=f_0$}}
\bezier{224}(63.00,30.00)(88.00,43.00)(113.00,30.00)
\put(88.00,33.00){\makebox(0,0)[cc]{$f_t$}}
\put(88.00,2.00){\makebox(0,0)[cc]{{\Large $2$-cell $f$}}}
\put(3.00,10.00){\line(1,0){40.00}}
\put(43.00,10.00){\line(0,1){40.00}}
\put(43.00,50.00){\line(-1,0){40.00}}
\put(3.00,50.00){\line(0,-1){40.00}}
\put(3.00,36.00){\vector(1,0){40.00}}
\end{picture}

\nopagebreak[4] Figure \theppp.
\end{center}
\addtocounter{ppp}{1}

By definition, a {\em positive path in a directed $2$-complex\/}
$\kk$ is a finite product of the form $p_1\cdots p_n$ ($n\ge1$),
where each factor $p_i$ ($1\le i\le n$) is either a $1$-path on
$\kk$, or a path of the form $f_t$ for some $f\in\F$, $t\in[0,1]$.
Of course, we assume that the terminal point of $p_i$ coincides
with the initial point of $p_{i+1}$ for any $1\le i<n$. The set
of all positive paths in $\kk$ defined in this way will be denoted
by $\Omega_+(\kk)$. Note that every $1$-path in $\kk$ is a positive
path in this sense.

Note that every positive path $p$ in $\Omega_+(\kk)$ can be uniquely
written in the {\em normal form\/}
$p=u_0f^{(1)}_{t_1}u_1\cdots f^{(n)}_{t_n}u_n$, where $u_i$ are $1$-paths
in $\kk$ ($0\le i\le n$), $f^{(i)}\in\F^-$, $t_i\in(0,1)$ ($1\le i\le n$).
To the path $p$ we assign a point $\xi(p)$ with coordinates
$(t_1,\ldots,t_n)$ in the $n$-cube $\ve(u_0)+f^{(1)}+\cdots+f^{(n)}+\ve(u_n)$
(see Figure \theppp).

\begin{center}
\unitlength=1.00mm
\special{em:linewidth 0.4pt}
\linethickness{0.4pt}
\begin{picture}(150.66,20.00)
\put(10.00,9.87){\circle*{1.33}}
\put(40.00,9.87){\circle*{1.33}}
\put(50.00,9.87){\circle*{1.33}}
\put(0.00,9.87){\circle*{1.33}}
\put(80.00,9.87){\circle*{1.33}}
\put(150.00,9.87){\circle*{1.33}}
\put(140.00,9.87){\circle*{1.33}}
\put(110.00,9.87){\circle*{1.33}}
\bezier{180}(10.00,10.00)(26.00,27.00)(40.00,10.00)
\bezier{180}(10.00,10.00)(26.00,-7.00)(40.00,10.00)
\bezier{192}(50.00,10.00)(65.00,29.00)(80.00,10.00)
\bezier{192}(50.00,10.00)(65.00,-9.00)(80.00,10.00)
\bezier{200}(110.00,10.00)(125.00,30.00)(140.00,10.00)
\bezier{196}(110.00,10.00)(127.00,-9.00)(140.00,10.00)
\linethickness{1.0pt}
\put(0.00,10.00){\line(1,0){10.00}}
\put(40.00,10.00){\line(1,0){10.00}}
\put(140.00,10.00){\line(1,0){10.00}}
\put(80.00,10.00){\line(1,0){8.00}}
\put(88.00,10.00){\line(0,0){0.00}}
\put(110.00,10.00){\line(-1,0){8.00}}
\linethickness{0.4pt}
\put(95.00,10.00){\makebox(0,0)[cc]{\dots}}
\linethickness{1.0pt}
\bezier{128}(10.00,10.00)(27.00,16.00)(40.00,10.00)
\bezier{136}(50.00,10.00)(64.00,2.00)(80.00,10.00)
\bezier{124}(110.00,10.00)(127.00,14.00)(140.00,10.00)
\linethickness{0.4pt}
\put(5.00,13.00){\makebox(0,0)[cc]{$u_0$}}
\put(45.00,13.00){\makebox(0,0)[cc]{$u_1$}}
\put(145.00,13.00){\makebox(0,0)[cc]{$u_n$}}
\put(25.00,9.00){\makebox(0,0)[cc]{$f^{(1)}_{t_1}$}}
\put(65.00,10.00){\makebox(0,0)[cc]{$f^{(2)}_{t_2}$}}
\put(125.00,9.00){\makebox(0,0)[cc]{$f^{(n)}_{t_n}$}}
\end{picture}

\nopagebreak[4] Figure \theppp.
\end{center}
\addtocounter{ppp}{1}

\begin{lm}
\label{OmSq}
The map $\xi$ is a bijection between $\Omega_+(\kk)$ and a realization of
$\Sq(\kk)$. For any positive path
\be{pth}
p=u_0f^{(1)}_{t_1}u_1\cdots f^{(n)}_{t_n}u_n
\ee
in $\Omega_+(\kk)$, where $u_i$ are $1$-paths in $\kk$ $(0\le i\le n)$ and
$t_i\in[0,1]$ $(1\le i\le n)$, the point $\xi(p)$ has coordinates
$(t_1,\ldots,t_n)$ in the $n$-cube $\ve(u_0)+f^{(1)}+\cdots+f^{(n)}+\ve(u_n)$.
\end{lm}

\proof
Any point $x$ in a realization of $\Sq(\kk)$ is an inner point of a unique cube
(the point of a cube of dimension $0$ is inner in itself by definition). Let $x$
be a point in the $n$-cube $\ve(u_0)+f^{(1)}+\cdots+f^{(n)}+\ve(u_n)$ with
coordinates $(t_1,\ldots,t_n)$, where $f^{(i)}\in\F^-$, $t_i\in(0,1)$ ($1\le i\le n$).
One can assign to it the positive path $p=u_0f^{(1)}_{t_1}u_1\cdots f^{(n)}_{t_n}u_n$
(written in the normal form). This defines a map $\zeta$ from a realization of
$\Sq(\kk)$ to $\Omega_+(\kk)$. It is obvious that $\xi$ and $\zeta$ are mutually
inverse. Thus $\xi$ is a bijection.

Now we have to check that if $p$ is written in the form (\ref{pth}), where
some subscripts $t_i$ are equal to $0$ or $1$, then the image of $p$ under
$\xi$ is defined by the same rule. Changing the reference point of a cube,
we may assume that $f^{(i)}\in\F^-$ for all $i$. Now one can pass from
(\ref{pth}) to the normal form of $p$ step by step replacing subpaths of the
form $f^{(i)}_0$ by $\topp{f^{(i)}}$ and $f^{(i)}_1$ by $\bott{f^{(i)}}$. Let
us analyze what happens at one elementary step. Without loss of generality
assume that a subpath of the form $f_0$ is replaced by $\topp{f}$. So our
path is written in two forms: $\cdots u_i\cdots$ and $\cdots u'f_0u''\cdots$,
where $u_i=u'\topp{f}u''$. The rule defining $\xi$ assigns two points to these
forms. One of them belongs to the cube $c'=\cdots+\ve(u_i)+\cdots$ and has
coordinates $(t_1,\ldots,t_n)$ in it. The other point belongs to the cube
$c=\cdots+\ve(u')+f+\ve(u'')+\cdots$ and has coordinates
$(t_1,\ldots,0,\ldots,t_n)$, where $0$ is inserted after $t_i$.
The cube $c'$ is the $(i+1)$th face of $c$. By definition of $\Sq(\kk)$, we
glue $c'$ and $\topp{c}_{i+1}$ isometrically. This means that these two
points in a realization of $\Sq(\kk)$ coincide.
\endproof

Lemma \ref{OmSq} shows that $\Omega_+(\kk)$ is a natural realization of the
semi-cubical complex $\Sq(\kk)$. Note that the empty paths correspond to
isolated points in the Squier complex.

Thus we have an equivalent definition of diagram groups of $\kk$ as the
fundamental groups of the space of positive paths in $\kk$. One possible way
of generalizing the definition of diagram groups and of defining ``continuous
versions" of the diagram groups would be to consider a more general spaces
than directed $2$-complexes, define ``positive paths" in a suitable way, and
then to consider the fundamental groups of the spaces of positive paths. They
may have certain properties in common with diagram groups.

Usually we shall not distinguish between the Squier complex $\Sq(\kk)$ and
its geometric realization.
\vspace{0.5ex}

We shall return to the Squier complexes and consider their universal covers
in Section \ref{ucd2c}.

\section{Theorems about isomorphism. The class of diagram groups}
\label{taicdg}

In this section, we show that diagram groups do not change much if we do
certain surgeries on directed $2$-complexes.

The following useful statement contains a directed $2$-complex analog of
Tietze transformations for group and semigroup presentations.

\begin{thm}
\label{tietze}
Let $\kk$ be a directed $2$-complex.

$1)$ Let $u$ be a non-empty $1$-path in $\kk$. Let $\kk'$ be the
directed $2$-complex obtained from $\kk$ by adding a new edge $e$
with $\iota(e)=\iota(u)$, $\tau(e)=\tau(u)$ and a new $2$-cell $f$
of the form $u=e$\ $($and also the inverse $2$-cell $f^{-1})$. Then
for every $1$-path $w$ in $\kk$, the diagram groups $\dd(\kk,w)$
and $\dd(\kk',w)$ are isomorphic.

$2)$ Suppose that $\kk$ is a union of two directed $2$-complexes
$\kk_1$ and $\kk_2$ such that all vertices and edges of $\kk$ are
both in $\kk_1$ and in $\kk_2$. Suppose that the top path
$\topp{f}$\ $($bottom path $\bott{f})$ of each positive $2$-cell
$f$ of $\kk_1$ is homotopic in $\kk_2$ to some path $u_f$\
$($resp., $v_f)$. Let us consider the directed complex $\kk'$ with
the same vertices and edges as $\kk$ and $2$-cells from $\kk_2$
together with all positive $2$-cells $u_f=v_f$ for all positive
$2$-cells $f$ from $\kk_1$\ $($plus the corresponding negative
cells $v_f=u_f)$. Then the diagram groups $\dd(\kk,w)$ and
$\dd(\kk',w)$ are isomorphic for every $1$-path $w$ in $\kk$.
\end{thm}

\proof
1. Indeed, there exists a natural embedding of $\dd(\kk,w)$ into
$\dd(\kk',w)$ which maps every reduced $(w,w)$-diagram over $\kk$ to
itself. In order to show that this map is surjective, notice that if
a $(w,w)$-diagram $\Delta$ over $\kk'$ does not contain edges labelled
by $e$ then it is a diagram over $\kk$. If $\Delta$ contains an edge
labelled by $e$ then this edge cannot be on $\topp{\Delta}$ or
$\bott{\Delta}$. Hence this edge is a common edge of the contours of
two cells $\pi$ and $\pi'$ in $\Delta$. Since $\kk'$ has only two
$2$-cells with $e$ on the boundary (namely, $f$ and $f^{-1}$), one of the
two cells $\pi$ or $\pi'$ is labelled by $f$ and another by $f^{-1}$ (the
edge labelled by $e$ is the top path of one of these cells and the bottom
path of another one). Hence $\pi$ and $\pi'$ form a dipole. This implies
that every reduced $(w,w)$-diagram over $\kk'$ contains no edges labelled
by $e$, and so it is a diagram over $\kk$. Hence by Theorem \ref{thnf}
the natural embedding of $\dd(\kk,w)$ into $\dd(\kk',w)$ is surjective.

2. By Theorem \ref{kash}, for any $2$-cell $f$ of $\kk_1$ there
exist diagrams $\Gamma_f$ and $\Delta_f$ over $\kk_2$ such that
the label of $\topp{\Gamma_f}$ is $\topp{f}$, the label of
$\bott{\Gamma_f}$ is $u_f$, the label of $\topp{\Delta_f}$ is
$\bott{f}$, the label of $\bott{\Delta_f}$ is $v_f$.

By $D_w$ (respectively, $D'_w$) we denote the set of all
$(w,w)$-diagrams over $\kk$ (respectively $\kk'$). We are going to
define two maps $\phi\colon D_w\to D'_w$, $\psi\colon D'_w\to
D_w$.

Let $\Xi\in D_w$. Let $\F_1^+$ be the set of positive $2$-cells in
$\kk_1$. For every cell $\pi$ in $\Xi$ with inner label $f\in
\F_1^+$, we do the following operation. First we cut $\pi$ into
three parts by connecting the initial vertex of $\pi$ with the
terminal vertex of $\pi$ by two simple curves, $p_1$ and $p_2$,
that have no intersections other than at the endpoints. We
enumerate the three parts from top to bottom and assume that $p_1$
is above $p_2$. Then we subdivide $p_1$ into edges and give them
labels such that $p_1$ will have label $u_f$. Similarly, we turn
$p_2$ into a path labelled by $v_f$. Now we insert the diagram
$\Gamma_f$ between the top path of $\pi$ and $p_1$. Analogously,
we insert $\Delta_f^{-1}$ between $p_2$ and the bottom path of
$\pi$. The space between $p_1$ and $p_2$ becomes a cell $u_f=v_f$,
which is a cell  $\kk'$. We can assign the inner label $f$ to it.
If the inner label of a cell $\pi$ of $\Xi$ is $f^{-1}\in\F_1^-$,
then we subdivide it in the same way to get the mirror image of
the diagram we had for cells with inner label $f$. (The inner
label for the cell in the middle will be $f^{-1}$.)

Every diagram $\Xi$ over $\pp$ now becomes a diagram over $\pp'$.
We denote it by $\phi(\Xi)$.

The map $\psi$ is defined similarly. Now if we have a diagram
$\Xi$ over $\kk'$, then we replace each of its cells $\pi$ of the
form $u_f=v_f$ ($f\in \F_1$) by the concatenation of three
diagrams. The first of them is $\Gamma_f^{-1}$, the third is
$\Delta_f$, and the second one is a cell with inner label $f$. We
do similar transformation with cells of the form $v_f=u_f$ whose
inner labels are negative.

The result of these replacements will be a diagram over $\kk$
denoted by $\psi(\Xi)$.

It follows from our construction that for any diagram $\Xi$ over
$\kk$, the diagram $\psi(\phi(\Xi))$ over $\kk$ is equivalent to
$\Xi$. This is so because after applying $\phi$ and then $\psi$ to
$\Xi$, we get a diagram with a number of subdiagrams of the form
$\Gamma^{\pm1}\Gamma^{\mp1}$ or $\Delta^{\pm1}\Delta^{\mp1}$.
Cancelling all the dipoles, we get the diagram $\Xi$ we had in the
beginning. Analogously, for any diagram $\Xi$ over $\kk'$, the
diagram $\phi(\psi(\Xi))$ over $\kk'$ will be equivalent to $\Xi$.
It is also clear that $\phi$ and $\psi$ preserve the operation of
concatenation of diagrams.

This means that maps $\phi$, $\psi$ induce homomorphisms of diagram
groups $\bar\phi\colon\dd(\kk,w)\to\dd(\kk',w)$ and
$\bar\psi\colon\dd(\kk',w)\to\dd(\kk,w)$. The fact about equivalence
of diagrams means that $\bar\phi$ and $\bar\psi$ are mutually inverse.
Thus they are isomorphisms and $\dd(\kk,w)\cong\dd(\kk',w)$.
\endproof

As an immediate application of Theorem \ref{tietze}, we obtain the
following statement about subdivisions of directed $2$-complexes.
Let $\kk$ be a directed $2$-complex and let $f$ be its $2$-cell.
Let us add a new edge $e$ to the complex with
$\iota(e)=\iota(\topp{f})=\iota(\bott{f})$ and
$\tau(e)=\tau(\topp{f})=\tau(\bott{f})$, remove the $2$-cells
$f^{\pm 1}$, and add new $2$-cells $f_1^{\pm1}$, $f_2^{\pm1}$,
where $f_1$, $f_2$ have the form $\topp{f}=e$ and $e=\bott{f}$,
respectively. This operation can be done for several positive
$2$-cells of $\kk$ at once. This simply means that we cut some
$2$-cells of $\kk$ into two parts. The resulting directed
$2$-complex $\kk'$ is called a {\em subdivision\/} of $\kk$.

\begin{lm}
\label{subd}
If $\kk$ is a directed $2$-complex, $w$ is a non-empty $1$-path in
$\kk$ and $\kk'$ is a subdivision of $\kk$, then the diagram groups
$\dd(\kk,w)$ and $\dd(\kk',w)$ are isomorphic.
\end{lm}

\proof
We use the notation from the paragraph preceding the formulation of
the lemma. Let $\kk''$ be the directed $2$-complex obtained from $\kk$
by adding the edge $e$ and the $2$-cells $f_1^{\pm1}$. By part 1) of
Theorem \ref{tietze}, $\dd(\kk,w)=\dd(\kk'',w)$. Now represent $\kk''$
as the union of $\kk_1$ and $\kk_2$, where $\kk_1$, $\kk_2$ have the
same vertices and edges, $\kk_1$ contains exactly two $2$-cells $f$,
$f^{-1}$, and $\kk_2$ contains all other $2$-cells of $\kk''$. Notice that
$\topp{f}$ is homotopic to $e$ in $\kk_2$ (because $\kk_2$ contains the
$2$-cell $f_1$). Hence by part 2) of Theorem \ref{tietze}, we can replace
$f^{\pm1}$ in $\kk''$ by $f_2^{\pm1}$, where $f_2$ has the form $e=\bott{f}$
without changing the diagram group with the base $w$. But the resulting
directed $2$-complex is precisely $\kk'$. Hence $\dd(\kk,w)\cong\dd(\kk',w)$.
\endproof

As an immediate corollary of Lemma \ref{subd} we get the following
statement.

\begin{thm}
\label{coinced}
The classes of diagram groups over semigroup presentations and diagram
groups of directed $2$-complexes coincide.
\end{thm}

\proof
Notice that complexes of the form $\kk_\pp$ corresponding to semigroup
presentations considered in \cite{GuSa97} are precisely the directed
$2$-complexes with one vertex in which different $2$-cells cannot have
the same top and bottom paths. We have already mentioned that we can only
consider directed $2$-complexes with one vertex (if we identify vertices we
preserve existing diagram groups but the set of diagram groups can
increase since the set of $1$-paths can increase). It is obvious that if
we subdivide each $2$-cell of $\kk$ twice (into three parts instead of two)
then we turn $\kk$ into a $\kk_\pp$ for some $\pp$. It remains to apply
Lemma \ref{subd}.
\endproof

\section{Morphisms of complexes and universal diagram groups}
\label{morph}

Let $\kk$, $\kk'$ be directed $2$-complexes. A {\em morphism\/} $\phi$
from $\kk$ to $\kk'$ is a map that takes vertices to vertices, edges
to non-empty $1$-paths and $2$-cells to $2$-paths and preserves the
functions $\iota$, $\tau$, $\topp{\cdot}$, $\bott{\cdot}$, $^{-1}$:

\begin{description}
\item[M1]\ For every edge $e$, $\phi(\iota(e))=\iota(\phi(e))$,
$\phi(\tau(e))=\tau(\phi(e))$,

\item[M2]\ For every $2$-cell $f$ of $\kk$,
$\phi(\topp{f})=\topp{\phi(f)}$, $\phi(\bott{f})=\bott{\phi(f)}$;
here for every $1$-path $p=e_1e_2\cdots e_k$ we set
$\phi(p)=\phi(e_1)\phi(e_2)\cdots\phi(e_n)$, where $e_i$ are edges
(the latter product exists because of M1\,).

\item[M3]\ For every $2$-cell $f$ of $\kk$, $\phi(f^{-1})=\phi(f)^{-1}$.
\end{description}

Every morphism $\phi\colon\kk\to\kk'$ induces a functor from the category
$\Pi(\kk)$ to $\Pi(\kk')$: the image of an atomic $2$-path $(p,f,q)$ is
$(\phi(p),\phi(f),\phi(q))$, where $(u,\delta_1\circ\cdots\circ\delta_n,v)$
for atomic $2$-paths $\delta_i=(p_i,f_i,q_i)$ is short for
$(up_1,f_1,q_1v)\circ\cdots\circ(up_n,f_n,q_nv)$.

A morphism $\phi\colon\kk\to\kk'$ also induces a functor from the diagram
groupoid $\dd(\kk)$ to $\dd(\kk')$. The image of a $(p,q)$-diagram $\Delta$
over $\kk$ is the $(\phi(p),\phi(q))$-diagram obtained from $\Delta$ by a)
replacing each edge that has label $e$ by a path labelled by $\phi(e)$, and
b) replacing each cell that has label $f$ by the diagram over $\kk'$
corresponding to the $2$-path $\phi(f)$. Both of these functors will be
denoted by $\phi$ as well. The restriction of $\phi$ onto a diagram group
$\dd(\kk,p)$ is a homomorphism that will be denoted by $\phi_p$.

For a non-empty $1$-path $p$ in $\kk$, we say that a morphism
$\phi\colon\kk\to\kk'$ is $p$-{\em nonsingular\/} if the induced
homomorphism $\phi_p$ is injective on $\dd(\kk,p)$. In that case
the subgroup $\phi_p(\dd(\kk,p))$ of $\dd(\kk')$ is called {\em
naturally embedded\/}. If the morphism $\phi$ is $p$-nonsingular
for every $p$, we call it {\em nonsingular\/}.

\begin{lm}
\label{emb1}
Let $\phi\colon\kk\to\kk'$ be a morphism of two directed $2$-comp\-lex\-es.

$1)$ If $\phi(\delta)$ is reduced for every reduced $2$-path $\delta$, and
$\phi(f)$ is not empty for every $2$-cell $f$ of $\kk$, then $\phi$ is
nonsingular.

$2)$ Suppose that $\phi$ is injective on the set of $2$-cells and $\phi(f)$
is a $2$-cell for every $2$-cell $f$ of $\kk$. Then $\phi$ is nonsingular.

$3)$ If a directed $2$-complex $\kk'$ is obtained by adding $2$-cells to a
directed $2$-complex $\kk$, then diagram groups of the form $\dd(\kk,p)$ are
naturally embedded into the diagram groups $\dd(\kk',p)$, for every $1$-path $p$.

$4)$ If a directed $2$-complex $\kk'$ is obtained by adding $2$-cells to a
directed $2$-complex $\kk$, and $\topp{f}$ is homotopic to $\bott{f}$ in
$\kk$ for every $2$-cell $f\in \kk'\setminus\kk$, then for every $1$-path
$p$ in $\kk$, the diagram group $\dd(\kk,p)$ is a retract of the diagram
group $\dd(\kk',p)$.
\end{lm}

\proof
1) Indeed, suppose that the kernel of $\phi_p$ is not trivial. Then it
contains a reduced $2$-path $\delta\ne p$ by Theorem \ref{thnf}. By the
assumption, $\phi(\delta)$ is reduced, and non-empty, so
$\phi_p(\delta)\ne1$ by Theorem \ref{thnf}, a contradiction.

2) It is easy to see that for every reduced $(p,p)$-diagram $\Delta$,
the diagram $\phi_p(\Delta)$ does not contain dipoles. It remains to
use part 1) of this lemma.

3) Immediately follows from 2).

4) The retraction $\psi$ is given as follows. Fix a
$(\topp{f},\bott{f})$-diagram $\Gamma_f$ over $\kk$ for every
$2$-cell $f\in\kk'\setminus\kk$ in such a way that inverse
diagrams correspond to inverse $2$-cells. Then for every
$(p,p)$-diagram $\Delta$ over $\kk'$, the diagram $\psi(\Delta)$
is obtained from $\Delta$ by inserting $\Gamma_f$ instead of every
cell in $\Delta$ labelled by $f\in\kk\setminus\kk'$. Clearly,
$\psi^2=\psi$.
\endproof

Let us call a directed $2$-complex $\kk$ {\em universal\/} if every
finite or countable directed $2$-complex maps nonsingularly into
$\kk$. A diagram group is called {\em universal\/} if it contains
copies of all countable diagram groups. A directed $2$-complex is
said to be {\em $2$-path connected} if all its non-empty $1$-paths
are homotopic to each other. By Theorem \ref{fg} it has  at most
one non-trivial diagram group up to isomorphism. Notice that if a
universal directed $2$-complex is $2$-path connected then its nontrivial
diagram group is universal because every at most countable diagram group
is (obviously) a diagram group of at most countable directed
$2$-complex.

\begin{lm}
\label{emb2}
Let $\uu$ be the directed $2$-complex
$$
\la x\mid x^m=x^n, x^m=x^n,\ldots\hbox{\ \rm for all\ }1\le m<n\ra
$$
$($every equality appears countably many times$)$. Then $\uu$ is universal.
\end{lm}

\proof
Let $\kk$ be at most countable directed $2$-complex and let $\kk_1$
be obtained from $\kk$ by identifying all its edges and vertices. By
$\phi$ we denote the natural morphism from $\kk$ to $\kk_1$. Then $\phi$
is nonsingular by Lemma \ref{emb1}, part 2. It is easy to see that $\uu$
can be obtained from $\kk_1$ by adding $2$-cells. Hence $\kk_1$ is a
subcomplex of $\uu$ and $\phi\colon\kk\to\uu$ is nonsingular by Lemma
\ref{emb1}, part 3.
\endproof

We can simplify the universal directed $2$-complex $\uu$ by using
part 2) of Theorem \ref{tietze}. For every $1\le n\le\infty$ let
$$
\hh_n=\la x\mid x^2=x,x=x,x=x,\ldots\ \ (n\hbox{ times})\ra.
$$

\begin{center}
\unitlength=1.00mm
\special{em:linewidth 0.4pt}
\linethickness{0.4pt}
\begin{picture}(29.00,35.00)
\put(1.00,19.00){\vector(1,1){14.00}}
\put(15.00,33.00){\vector(1,-1){14.00}}
\put(1.00,19.00){\vector(1,0){28.00}}
\bezier{144}(1.00,19.00)(15.00,7.33)(29.00,19.00)
\bezier{216}(1.00,19.00)(15.00,-4.33)(29.00,19.00)
\put(15.00,4.50){\makebox(0,0)[cc]{$\ldots$}}
\bezier{300}(1.00,19.00)(15.00,-15.67)(29.00,19.00)
\bezier{120}(1.00,19.00)(13.99,13.95)(29.00,19.00)
\put(16.00,16.50){\vector(1,0){0.00}}
\put(16.00,13.15){\vector(1,0){0.00}}
\put(16.00,7.33){\vector(1,0){0.00}}
\put(16.00,1.67){\vector(1,0){0.00}}
\end{picture}

\nopagebreak[4] Figure \theppp.
\end{center}

\noindent
This complex contains one vertex, one edge, and $n+1$ positive
$2$-cells, one of which is the Dunce hat and all others are spheres.
It is obtained from the plane diagram on Figure \theppp\ by
identifying all edges according to their directions. For example, the
complex $\hh_0$ is the Dunce hat (see Figure 1).

\addtocounter{ppp}{1}

\begin{lm}
\label{emb3}
The directed $2$-complex $\hh_\infty$ is universal. In particular, every
countable diagram group embeds into the diagram group $\dd(\hh_\infty,x)$.
\end{lm}

\proof In fact we shall show that $\dd(\uu,x)$ is isomorphic to
$\dd(\hh_\infty,x)$. Let us denote by $\qq$ the complex obtained
from $\uu$ by removing the cell $x^2=x$ and its inverse. Then
$\hh_0\cup\qq=\uu$. It is easy to see that every non-empty $1$-path
in $\hh_0$ is homotopic to $x$. Therefore, let us replace each
$2$-cell $x^m=x^n$ in $\uu$ by a cell $x=x$, and obtain a complex
$\hh_\infty$. By part 2) of Theorem \ref{tietze}, the diagram
groups of $\uu$ and $\hh_\infty$ are isomorphic. The proof of part
2) of Theorem \ref{tietze} actually gives us a nonsingular
morphism from $\uu$ into $\hh_\infty$. It remains to use Lemma
\ref{emb2}.
\endproof

The diagram group $\dd(\hh_\infty,x)$ is not even finitely
generated (see Example \ref{ex33} below).

Our next goal is to show that already the group $\dd(\hh_1,x)$ is
universal. This will follow from Lemma \ref{emb3} and the fact
that there exists a nonsingular morphism from $\hh_\infty$ to
$\hh_1$. The group $\dd(\hh_1,x)$ is finitely presented (see
Example \ref{ex33}) and has a nice structure (see \cite{GuSa02b}).

\begin{lm}
\label{emb4}
There exists a nonsingular morphism from $\hh_\infty$ to $\hh_2$.
\end{lm}

\proof
Let us label the positive $2$-cells of $\hh_\infty$ by $f_0$, $f_1$,
\dots\,, where $f_0$ is the cell $x^2=x$. The complex $\hh_2$ has
positive $2$-cells $f_0$, $f_1$, $f_2$. Consider the morphism $\phi$
from $\hh_\infty$ to $\hh_2$ which takes the edge $x$ to $x$, $f_0$
to $f_0$, and each $f_i$, $i\ge 1$, to the $2$-path
\be{67}
(1,f_1,1)^i\circ (1,f_2,1)\circ (1,f_1,1)^i.
\ee
We need to show that $\phi_x$ is injective (then every $\phi_{x^k}$ will
be injective too because all local groups in the diagram groupoid
$\dd(\hh_\infty)$ are conjugate by Corollary \ref{cy1}).

Indeed, let $\Delta$ be a nontrivial reduced diagram over
$\hh_\infty$. Then $\bar\Delta=\phi(\Delta)$ is obtained by
replacing every cell with label $f_i$ by the diagram $\Gamma_i$
corresponding to the $2$-path (\ref{67}) (see Figure \theppp) and
each cell with label $f_i^{-1}$ by the diagram $\Gamma_i^{-1}$.

\begin{center}
\unitlength=0.5mm
\special{em:linewidth 0.4pt}
\linethickness{0.4pt}
\begin{picture}(49.00,90.00)
\bezier{208}(-15.00,45.00)(27.00,60.33)(69.00,45.00)
\bezier{204}(-15.00,45.00)(27.00,31.67)(69.00,45.00)
\bezier{348}(-15.00,45.00)(27.00,83.33)(69.00,45.00)
\bezier{552}(-15.00,45.00)(27.00,110.67)(69.00,45.00)
\bezier{736}(-15.00,45.00)(27.00,134.67)(69.00,45.00)
\put(27.00,83.00){\makebox(0,0)[cc]{$\ldots$}}
\put(27.00,58.00){\makebox(0,0)[cc]{$f_1$}}
\put(27.00,71.00){\makebox(0,0)[cc]{$f_1$}}
\put(27.00,45.00){\makebox(0,0)[cc]{$f_2$}}
\bezier{348}(-15.00,45.00)(27.00,7.67)(69.00,45.00)
\bezier{496}(-15.00,45.00)(27.00,-12.67)(69.00,45.00)
\bezier{712}(-15.00,45.00)(27.00,-40.67)(69.00,45.00)
\put(27.00,10.00){\makebox(0,0)[cc]{$\ldots$}}
\put(27.00,21.00){\makebox(0,0)[cc]{$f_1$}}
\put(27.00,33.00){\makebox(0,0)[cc]{$f_1$}}
\end{picture}

\nopagebreak[4] Figure \theppp.
\end{center}
\addtocounter{ppp}{1}

It is sufficient to show that $\bar\Delta$ is also nontrivial. Any
diagram over $\hh_n$ (for any $n$) can be uniquely decomposed into
subdiagrams of the following two types. Each subdiagram of the
first type is an $(x^2,x)$-cell or its mirror image. Each subdiagram
of the second type is a maximal $(x,x)$-subdiagram, which is a product
of $(x,x)$-cells only.

If we decompose $\bar\Delta$ in such a way, then we see that each
subdiagram of the second type is the image of a subdiagram of the
second type in $\Delta$ under the mapping $\phi$. Since $\Delta$
has no dipoles, each of its subdiagram of the second type is a
product of $(x,x)$-cells, where the word formed by their labels is
a freely irreducible group word $v$ over the alphabet
$\{f_i^{\pm1}\mid i\ge1\}$. Notice that the $\phi$-image of this
subdiagram is a product of $(x,x)$-cells such that the word formed by
their inner labels is $\bar v$, where $\bar f_i=f_1^if_2f_1^i$.

Since the map $f_i\mapsto\bar f_i$ is an embedding of the free
group with generators $f_i$ ($i\ge1$) into the free group with
two generators $f_1$, $f_2$, the word $\bar v$ will be non-empty
after all free cancellations are made in it. Thus if we reduce
dipoles in each subdiagram of the second type in $\bar\Delta$, the
resulting diagram $\Delta'$ will be reduced. Indeed, we have
cancelled all the dipoles formed by $(x,x)$-cells. No dipoles
between cells that correspond to $x^2=x$ may appear because each
subdiagram of the second type in $\bar\Delta$ remains nontrivial
after all cancellations. The diagram $\Delta'$ contains the same
number of cells labelled by $f_0^{\pm 1}$ as $\Delta$ and at least
as many $(x,x)$-cells as $\Delta$. Hence $\Delta'$ is nontrivial.
\endproof

\begin{lm}
\label{g2tog1}
There exists a nonsingular morphism from $\hh_2$ into $\hh_1$.
\end{lm}

\proof
The idea is similar to the one of the previous lemma. We keep notation
for $2$-cells of $\hh_2$ from that lemma. Letters $f_i$ ($i=1,2$) will
be also used to denote the atomic $2$-paths $(1,f_i,1)$ and the
corresponding $(x,x)$-diagrams that consist of one cell labelled by $f_i$.
By $a$ we denote any nontrivial reduced $(x,x)$-diagram over
$\hh_0=\la x\mid x^2=x\ra$ and one of the corresponding $2$-paths in
$\hh_0$. Any two $(x,x)$-diagrams can be concatenated. So each word in
$f_1^{\pm1}$, $f_2^{\pm1}$, $a^{\pm1}$ denotes some $(x,x)$-diagram.

Now we use the morphism $\psi\colon\hh_2\to\hh_1$ that takes the
edge $x$ to $x$, $f_1$ to $f_1af_1$, $f_2$ to $f_1^2af_1^2$.

Let $\Delta$ be a nontrivial reduced diagram over $\hh_2$. Let
$\tilde\Delta$ be the diagram obtained from $\psi(\Delta)$ by cancelling
all dipoles of $(x,x)$-cells. As in Lemma \ref{emb4}, we subdivide $\Delta$
into subdiagrams of the two types. Let $\Gamma$ be a subdiagram of the second
type. It is a product of cells with inner labels $f_1^{\pm1}$, $f_2^{\pm1}$.
Let $v$ be the word that is the product of these labels. Clearly, this
is a freely reduced word in $f_1^{\pm1}$, $f_2^{\pm1}$. After we replace
$v$ by $\bar v$, where $\bar f_1=f_1af_1$, $\bar f_2=f_1^2af_1^2$ and then
freely reduce the result, we get a word of the form
\be{frm}
f_1^{s_0}a^{k_1}f_1^{s_1}\cdots a^{k_r}f_1^{s_r},
\ee
where $r$ is the length of $v$. Note that $s_0\ne0$, $s_r\ne0$,
$k_1,\ldots,k_r=\pm1$. Note also that none of the occurrences of letters
$a^{\pm 1}$ in $\bar v$ disappears after the reduction, hence only
occurrences of the letter $f_1^{\pm1}$ can disappear. Thus the
$\psi$-image of any subdiagram of the second type after cancelling all
dipoles of $(x,x)$-cells becomes reduced and of the form (\ref{frm}) as well.
After we cancel all $(x,x)$-dipoles in the subdiagrams $\psi(\Gamma)$ for
all maximal subdiagrams $\Gamma$ of $\Delta$ of the second type, we would
not have any more $(x,x)$-dipoles. Hence we shall get the diagram
$\tilde\Delta$.

Since $s_0$ and $s_r$ are always non-zero, the $2$-cells labelled by
$f_0^{\pm 1}$ cannot form a dipole in $\tilde\Delta$. Therefore,
$\tilde\Delta$ is reduced. Since the number of cells in
$\tilde\Delta$ is at least the same as in $\Delta$, the diagram
$\tilde\Delta$ is nontrivial.
\endproof

\begin{thm}
\label{th565}
The directed $2$-complex $\hh_1$ is universal. Hence the group
${\cal G}_1=\dd(\hh_1,x)$ contains copies of all countable diagram
groups. This group is finitely presented and even of type $\fff_\infty$.
\end{thm}

\proof
The first statement follows from Lemmas \ref{emb3}, \ref{emb4}, \ref{g2tog1}.
The fact that the group ${\cal G}_1$ is of type $\fff_\infty$ follows from
Theorem \ref{KGfinite} below. It can also be deduced from results of
\cite{Far00}. In Example \ref{ex33} below, we shall compute a presentation
of ${\cal G}_1$. It has $6$ generators and $18$ defining relations.
\endproof

Theorem \ref{th565} gives an example of a universal directed $2$-complex
with one edge and two positive $2$-cells. The following theorem shows that
a complex with one edge and one positive $2$-cell can be universal as well.

\begin{thm}
\label{th24}
The directed $2$-complex $\fft=\la y\mid y^3=y^2\ra$ is universal
$($this complex is obtained from Figure {\rm\theppp}\ by identifying
all edges according to their directions$)$. Its diagram group
$\dd(\fft,y^2)$ is also a universal group of type $\fff_\infty$. It has
the following Thompson-like presentation:
$$
\la x_0, x_1,\dots, y_0, y_1,\dots\mid x_j^{x_i}=x_{j+1},y_j^{x_i}=y_{i+1},
\ 0\le i<j-1\ra.
$$
\end{thm}

\begin{center}
\unitlength=1.00mm
\special{em:linewidth 0.4pt}
\linethickness{0.4pt}
\begin{picture}(62.00,24.00)
\put(1.00,2.00){\circle*{1.33}}
\put(21.00,22.00){\circle*{1.33}}
\put(41.00,22.00){\circle*{1.33}}
\put(61.00,2.00){\circle*{1.33}}
\put(31.00,2.00){\circle*{1.33}}
\put(1.00,2.00){\vector(1,1){10.00}}
\put(11.00,12.00){\line(1,1){10.00}}
\put(21.00,22.00){\line(1,0){20.00}}
\put(41.00,22.00){\vector(1,-1){10.00}}
\put(51.00,12.00){\line(1,-1){10.00}}
\put(1.00,2.00){\line(1,0){30.00}}
\put(31.00,2.00){\line(1,0){30.00}}
\put(31.00,22.00){\vector(1,0){0.00}}
\put(16.00,2.00){\vector(1,0){0.00}}
\put(46.00,2.00){\vector(1,0){0.00}}
\end{picture}

\nopagebreak[4] Figure \theppp.
\end{center}
\addtocounter{ppp}{1}

\proof
Let us consider the two diagrams $A$ and $B$ over $\fft$ on Figure \theppp.
Here $A$ is a $(y^8,y^4)$-diagram and $B$ is a $(y^4,y^4)$-diagram. These
diagrams correspond to the following $2$-paths on $\fft$. Let $p_{i,j}$
($i,j\ge0$) denote the atomic $2$-path $(y^i,y^3=y^2,y^j)$. Then $A$
corresponds to the $2$-path $\alpha=p_{33}^{-1}p_{15}p_{41}p_{04}p_{30}p_{11}$,
whereas $B$ corresponds to $\beta=p_{01}p_{10}^{-1}$. Let us consider the
morphism $\gamma$ from $\hh_1$ to $\fft$ which takes the edge $x$ to the
$1$-path $y^4$, the $2$-cell $f_0$ to $\alpha$, and the $2$-cell $f_1$ to
$\beta$. We are going to show that $\gamma$ is nonsingular. As before, it is
enough to show that it is $x$-nonsingular.

Let $\Delta$ be any reduced $(x,x)$-diagram over $\hh_1$. The diagram
$\gamma(\Delta)$ is obtained as follows. First we subdivide each edge
labelled by $x$ into $4$ parts and label each of them by $y$. Then each
$(x^2,x)$-cell becomes a $(y^8,y^4)$-cell. Every $(x^2,x)$-cell with inner
label $f_0$ is replaced by $A$. A mirror image of such a cell is replaced
by $A^{-1}$. Similarly, any $(x,x)$-cell of $\Delta$ becomes a $(y^4,y^4)$-cell,
so we replace by $B^{\pm1}$ all $(x,x)$-cells labelled by $f_1^{\pm1}$. After
all these replacements, we get a $(y^4,y^4)$-diagram $\hat\Delta$ over $\fft$.

We shall show that $\hat\Delta$ is reduced, and then apply Lemma \ref{emb1},
part 1). Note that each of the diagrams $A$, $B$ has no dipoles so a dipole in
$\hat\Delta$, if it occurs, must belong to different subdiagrams of the form
$A^{\pm1}$, $B^{\pm1}$. Suppose that the upper cell of the dipole is contained in
$A^{\pm1}$. From the structure of $A$ it is obvious that this cell must be a
$(y^3,y^2)$-cell. So it cannot form a dipole with a cell from $B^{\pm1}$. The
lower cell of the dipole is a $(y^2,y^3)$-cell.

\begin{center}
\unitlength=1.00mm
\special{em:linewidth 0.4pt}
\linethickness{0.4pt}
\begin{picture}(126.33,23.67)
\put(80.33,12.00){\circle*{1.33}}
\put(95.33,12.00){\circle*{1.33}}
\put(110.33,12.00){\circle*{1.33}}
\put(125.33,12.00){\circle*{1.33}}
\bezier{200}(80.33,12.00)(95.33,32.00)(110.33,12.00)
\bezier{212}(95.33,12.00)(110.33,-10.00)(125.33,12.00)
\put(88.33,20.00){\circle*{1.33}}
\put(102.33,20.00){\circle*{1.33}}
\put(102.33,4.00){\circle*{1.33}}
\put(118.33,4.00){\circle*{1.33}}
\put(4.33,1.67){\circle*{1.33}}
\put(19.33,1.67){\circle*{1.33}}
\put(34.33,1.67){\circle*{1.33}}
\put(49.33,1.67){\circle*{1.33}}
\put(64.33,1.67){\circle*{1.33}}
\bezier{180}(19.33,1.67)(34.33,18.67)(49.33,1.67)
\put(27.33,8.00){\circle*{1.33}}
\put(41.33,8.67){\circle*{1.33}}
\bezier{152}(27.33,8.00)(10.33,18.67)(4.33,1.67)
\bezier{152}(41.33,8.67)(60.33,17.67)(64.33,1.67)
\put(7.33,7.67){\circle*{1.33}}
\put(62.33,6.67){\circle*{1.33}}
\put(18.33,11.67){\circle*{1.33}}
\put(51.33,11.67){\circle*{1.33}}
\bezier{228}(7.33,7.67)(16.33,34.67)(27.33,8.67)
\bezier{240}(41.33,8.67)(54.33,35.67)(62.33,6.67)
\put(11.50,17.00){\circle*{1.33}}
\put(22.33,17.67){\circle*{1.33}}
\put(46.33,17.67){\circle*{1.33}}
\put(58.33,17.67){\circle*{1.33}}
\put(34.33,17.67){\circle*{1.33}}
\put(4.33,1.67){\line(1,0){60.33}}
\put(22.33,17.67){\line(1,0){24.00}}
\put(80.33,12.00){\line(1,0){45.00}}
\end{picture}

\nopagebreak[4] Figure \theppp.
\end{center}
\addtocounter{ppp}{1}

There are two types of vertices in $\hat\Delta$. Vertices of the
first type (we call them {\em red\/}) are the images of vertices
of $\Delta$. The other vertices are called {\em green\/}. It is
easy to see that $B$ has only two red vertices, $\iota(B)$ and
$\tau(B)$. The diagram $A$ has exactly three red vertices:
$\iota(A)$, $\tau(A)$, and the middle point of the top path of
$A$. The middle point of the bottom path of $A$ is green. So we
have to consider two cases for the two cells that form the dipole.
In the first case the middle point of the common part of the
boundary of the cells forming a dipole is red, in the second case
it is green.

In the first case the upper cell of the dipole is contained in a
copy of $A^{-1}$ and the lower cell is contained in a copy of $A$.
Denote these copies by $\Gamma_1$, $\Gamma_2$, respectively. We
claim that the bottom path of $\Gamma_1$ coincides with the top
path of $\Gamma_2$. Indeed, $\Gamma_1$ is the $\gamma$-image
of an $(x,x^2)$-cell $\pi_1$ in $\Delta$ and $\Gamma_2$ is the
$\gamma$ image of an $(x^2,x)$-cell $\pi_2$ in $\Delta$. The
bottom path of $\pi_1$ was subdivided into $8$ parts. The same is
true for the top path of $\pi_2$. The product of the $4$th and the
$5$th of these parts is the same for both $\pi_1$ and $\pi_2$
since it is the common boundary of the cells forming the dipole.
This can happen only if the bottom path of $\pi_1$ coincides with
the top path of $\pi_2$. But in this case we have a dipole in
$\Delta$. This contradicts the assumption that $\Delta$ is
reduced.

In the second case the upper cell of the dipole is contained in a
copy of $A$ and the lower cell is contained in a copy of $A^{-1}$.
We also denote these copies by $\Gamma_1$, $\Gamma_2$,
respectively. The bottom path of $\Gamma_1$ is the image of an edge
in $\Delta$. This edge is subdivided into $4$ parts. The product of
its second and third part is the common boundary of the cells of the
dipole. The same is true for the top path of $\Gamma_2$. So the bottom
of $\Gamma_1$ coincides with the top of $\Gamma_2$ since they must be
images of the same edge in $\Delta$. In this case an $(x^2,x)$-cell
forms a dipole in $\Delta$ with an $(x,x^2)$-cell, a contradiction.

If the lower cell of the dipole is contained in $A^{\pm1}$, then
the same arguments are applied. So to finish the proof, let us
assume that the dipole in $\hat\Delta$ is formed by two cells that
are contained in copies of $B^{\pm1}$. Suppose that the upper cell
of the dipole belongs to a copy of $B$. Thus it is a
$(y^2,y^3)$-cell. Note that the leftmost point of it is green and
the rightmost point is red. The lower cell must be a
$(y^3,y^2)$-cell with the corresponding points of the same colour.
Thus the lower cell is contained in a copy of $B^{-1}$. As in the
previous paragraph, we see from this fact that the last $3$ of $4$
sections of some edges in $\Delta$ coincide. Then these edges also
coincide and so $\Delta$ has a dipole that consists of two
$(x,x)$-cells. The case when the upper cell of the dipole belongs
to a copy of $B^{-1}$ is quite analogous. Thus $\gamma$ is
nonsingular.

That diagram groups of $\fft$ are of type $\fff_\infty$ follows
directly from \cite{Far00}. The presentation of $\dd(\fft,y^2)$ is
found in \cite[page 114]{GuSa97}. It can be found using Theorem
\ref{th11} below.
\endproof

In the next section we shall construct a universal directed $2$-complex
whose diagram groups have very simple presentations.

\section{Presentations of diagram groups}
\label{complete}

In \cite[Section 9]{GuSa97}, we showed how to find nice
presentations of diagram groups of the so called complete string
rewriting systems. Here we shall generalize these results for
diagram groups of directed $2$-complexes.

We start with a definition of a complete directed $2$-complex.
Throughout this section, $\kk$ is a directed $2$-complex with the
set of edges $\E$, set of $2$-cells $\F$ and a fixed set of positive
$2$-cells $\F^+$.

Let $p$, $q$ be $1$-paths in $\kk$. We write $p\tos q$ if $p\ne q$ and
there exists a positive $2$-path $\delta$ with $\topp{\delta}=p$ and
$\bott{\delta}=q$.

We say that $\kk$ is {\em Noetherian\/} if every sequence of $1$-paths
$p_1\tos p_2\tos\cdots$ terminates.

We say that $\kk$ is {\em confluent\/}, if for every two positive
$2$-paths $\delta_1$, $\delta_2$ with $\topp{\delta_1}=\topp{\delta_2}$
there exist two positive $2$-paths $\delta_1\circ\delta_1'$ and
$\delta_2\circ\delta_2'$ such that $\bott{\delta_1'}=\bott{\delta_2'}$.
In that case we say that $\delta_1$ and $\delta_2$ can be {\em extended
to a diamond\/}.

If a directed $2$-complex is Noetherian and confluent, then we say that
$\kk$ is {\em complete\/}.

It is easy to see that if $\kk=\kk_\pp$ for some complete string
rewriting system $\pp$, then $\kk$ is complete.

Let $\delta_1$, $\delta_2$ be two positive atomic $2$-paths on
$\kk$. Assume that $\topp{\delta_i}\ne\bott{\delta_i}$ for each
$i=1,2$. Suppose that one of the two cases hold:
\begin{enumerate}
\item $\delta_1=(1,f_1,q)$, $\delta_2=(p,f_2,1)$, where
$\topp{f_1}=ps$, $\topp{f_2}=sq$ for some non-empty $1$-path $s$;
\item $\topp{f_2}$ is a subpath of $\topp{f_1}$ and $f_1\ne f_2$.
\end{enumerate}
Then we say that $\delta_1$ and $\delta_2$ form a {\em critical
pair\/}. The diagrams representing these cases are shown on Figure
\theppp.

\begin{center}
\unitlength=1.00mm
\special{em:linewidth 0.4pt}
\linethickness{0.4pt}
\begin{picture}(125.00,27.33)
\put(5.67,14.67){\line(1,0){53.67}}
\bezier{200}(5.67,14.67)(21.33,34.00)(37.33,14.67)
\bezier{212}(27.00,14.67)(46.33,-6.33)(59.33,14.67)
\put(21.33,19.00){\makebox(0,0)[cc]{$f_1^{-1}$}}
\put(44.00,8.67){\makebox(0,0)[cc]{$f_2$}}
\put(21.00,12.00){\makebox(0,0)[cc]{$p$}}
\put(34.00,13.00){\makebox(0,0)[cc]{$s$}}
\put(50.00,12.00){\makebox(0,0)[cc]{$q$}}
\put(78.33,14.67){\line(1,0){44.33}}
\bezier{268}(78.67,14.67)(99.67,40.33)(122.33,14.67)
\bezier{164}(88.67,14.67)(103.00,-2.00)(113.00,14.67)
\put(99.67,20.67){\makebox(0,0)[cc]{$f_1^{-1}$}}
\put(101.33,11.00){\makebox(0,0)[cc]{$f_2$}}
\put(38.67,16.67){\vector(0,1){5.67}}
\put(41.00,19.00){\makebox(0,0)[cc]{$\delta_1$}}
\put(15.00,12.67){\vector(0,-1){5.33}}
\put(12.33,10.00){\makebox(0,0)[cc]{$\delta_2$}}
\put(122.67,15.67){\vector(0,1){5.67}}
\put(125.00,18.00){\makebox(0,0)[cc]{$\delta_1$}}
\put(86.00,13.33){\vector(0,-1){5.33}}
\put(83.33,10.67){\makebox(0,0)[cc]{$\delta_2$}}
\end{picture}

\nopagebreak[4] Figure \theppp.
\end{center}
\addtocounter{ppp}{1}

We say that the critical pair can be {\em resolved\/} if it can be
extended to a diamond.

For string rewriting systems, it is known (Newman's lemma \cite{GuSa97})
that a Noetherian string rewriting system is complete if and only if
every critical pair can be resolved. One can similarly prove that a
Noetherian directed $2$-complex is complete if and only if every
critical pair of its positive atomic $2$-paths can be resolved.

A $1$-path $p$ in a complete directed $2$-complex $\kk$ is called
{\em irreducible\/} if $p\tos q$ is impossible (that is, $p$ cannot be
changed by any positive $2$-path). It is easy to see that every
$1$-path $p$ in a complete directed $2$-complex $\kk$ is homotopic
to a unique irreducible $1$-path $\bar p$, which is called the
{\em irreducible form\/} of $p$.

From Lemma \ref{emb1}, part 4), one can almost immediately deduce
the following important statement.

\begin{lm}
\label{retract}
Every diagram group of a directed $2$-complex $\kk$ is a retract of a
diagram group of a complete directed $2$-complex $\kk'\supseteq\kk$.
The number of classes of homotopic $1$-paths in $\kk'$ is the same as
in $\kk$. If $\kk$ is finite and it has finitely many classes of homotopic
$1$-paths, then $\kk'$ is also finite.
\end{lm}

\proof
Let $G=\dd(\kk,p)$. Let us fix some total well ordering on the set of
edges of $\kk$. Then we can introduce the ShortLex order on $1$-paths
of $\kk$.

Let us change the orientation on the set of $2$-cells of $\kk$ as follows.
For every positive cell $f\in \kk$, if $\topp{f}$ is smaller than $\bott{f}$
in ShortLex, we call $f$ negative and $f^{-1}$ positive. This operation does
not change the diagram groups of the $2$-complex and the classes of homotopic
$1$-paths (see Remark \ref{rk33}).

In every class $W$ of homotopic $1$-paths of $\kk$ choose the
ShortLex smallest $1$-path $p(W)$. Now add cells to $\kk$ as
follows. First for every edge $e$ in $\kk$, we add a positive
$2$-cell $f_e$ of the form $e=p(W)$, where $W$ is the class of
homotopic $1$-paths containing $e$. We also add the inverse of that
$2$-cell. Now let $U$ and $V$ be two classes of homotopic $1$-paths
in $\kk$ such that the product $p(U)p(V)$ exists (that is,
$\tau(p(U))=\iota(p(V))$) and let $W$ be the class of homotopic
$1$-paths that contains $p(U)p(V)$. Add a positive cell $f_{U,V}$
with top path $p(U)p(V)$ and bottom path $p(W)$ (also add the
corresponding negative $2$-cell). As a result of these operations,
the number of classes of homotopic $1$-paths does not change.

The resulting complex $\kk'$ is clearly Noetherian. Indeed, for
every positive $2$-cell $f$ of $\kk'$, $\topp{f}\ge\bott{f}$ in
the ShortLex order.

Every $1$-path of a class $W$ is connected to $p(W)$ by a positive
$2$-path, which consists of atomic $2$-paths corresponding to the cells of
the form $f_e$ and $f_{U,V}$. (This can be easily proved by
induction on the length of the $1$-path.) Hence the complex $\kk'$
is confluent. Indeed, for every two positive atomic $2$-paths
$\delta_1$,  $\delta_2$ in $\kk$ with $\topp{\delta_1}=\topp{\delta_2}$,
their bottom $1$-paths are homotopic. So one can reduce each of
them to the same $1$-path and complete the diamond. Thus $\kk'$ is a
complete directed $2$-complex.

The complex $\kk'$ is also confluent. Indeed, for every two positive
atomic $2$-paths $\delta_1$, $\delta_2$ in $\kk$ with
$\topp{\delta_1}=\topp{\delta_2}$, the $1$-paths $\bott{\delta_1}$
and $\bott{\delta_2}$ are in the same class $W$ of homotopic $1$-paths.
Now using the new cells $f_e$ and $f_{U,V}$, one can reduce each of these
$1$-paths to $p(W)$ and complete the diamond. Thus $\kk'$ is a complete
directed $2$-complex.

By Lemma \ref{emb1}, part 4), $G$ is a retract of $\dd(\kk',p)$.
The last statement of the theorem obviously holds because we add
only finitely many cells.
\endproof

Since we are looking for nice presentations of diagram groups,
which are fundamental groups of Squier complexes, it is natural to
start with finding nice spanning forests in $\Sq(\kk)$. It can be
done in the case when $\kk$ is complete (and in some other cases
which we do not discuss here). Recall that a {\em spanning forest\/}
of a $1$-complex $\sss$ is a forest whose intersection with every
connected component of $\sss$ is a spanning tree in that component.

\begin{df}
\label{dfm}
{\rm
Let $\kk$ be a complete directed $2$-complex. A spanning forest $T$ in
$\Sq(\kk)$ is called a {\em left forest\/} whenever the following two
conditions hold:
\begin{description}
\item[F1]\ for any edge $e=(p,f,q)$ in $T$, the $1$-path $p$ is
irreducible;
\item[F2]\ if an edge $e=(p,f,q)$ belongs to $T$, then any
edge of the form $(p,f,q')$ also belongs to $T$.
\end{description}
Analogously one can define a {\em right forest\/}.
}
\end{df}

Because of the property F2, we will often use the notation $(u,f,*)$ when
we mention an edge of a left forest. Analogously, $(*,f,v)$ will be used
for edges from a right forest.

\begin{lm}
\label{comp-for}
If $\kk$ is a complete directed $2$-complex, then $\Sq(\kk)$ has a left
forest and a right forest.
\end{lm}

\proof
In each connected component of $\Sq(\kk)$ we choose the vertex which is an
irreducible $1$-path. If $p$ is not irreducible, then we find its shortest
initial segment $p'$, which is not irreducible. Let $p=p'v$ for some $v$.
By definition, $p'$ can be reduced so it has a subpath of the form $\bott{f}$
for some negative cell $f$ of $\kk$, where $\bott{f}\neq\topp{f}$ (there are
possibly many ways to choose $f$ with the above properties but we choose
{\bf one} of them arbitrarily). Obviously, this subpath is a suffix of $p'$.
Hence $p'=u\bott{f}$, where $u$ must be irreducible. Thus to every vertex
$p$ of $\Sq(\kk)$ that is not reducible, we can assign an edge $e=(u,f,v)$.
Let us consider the subgraph $T$ of the $1$-skeleton of $\Sq(\kk)$ that
contains all vertices and all the edges of the form $e^{\pm1}$, where $e$
was assigned to some $p$. We leave it as an exercise for the reader to check
that $T$ is a spanning forest and that it satisfies conditions F1, F2 (see
the proof of \cite[Lemma 9.4]{GuSa97}). A right forest in $\Sq(\kk)$ is
constructed in a similar way.
\endproof

\begin{rk}
\label{props}
{\rm
Let $\kk$ be a complete directed $2$-complex. In general, the way to
construct a left (right) forest from the proof of Lemma \ref{comp-for}
is not unique. However, if the second case of the critical pair from its
definition never occurs (that will be the case in all the examples
considered below), then it is not difficult to prove that the left forest
in $\Sq(\kk)$ is unique and consists of all edges $(u,f,v)^{\pm 1}$,
where $f\in \F^-$, $\bott{f}\ne\topp{f}$, and every proper initial subpath
of $u\bott{f}$ is irreducible.
}
\end{rk}

Let us fix a left forest $T_l$ and a right forest $T_r$ in $\Sq(\kk)$.
Then for every vertex $p$ in $\Sq(\kk)$, where $p\ne\bar p$, there exists
a unique negative edge $e\in T_l$ (resp., $e\in T_r$) going into $p$.
Indeed, otherwise there would be two different paths in $T_l$ (resp., $T_r$)
that consist of positive edges and connect $p$ with $\bar p$. We shall say
that $e$ is {\em assigned\/} to $p$.

The following theorem is a translation of \cite[Theorem 9.5]{GuSa97}.
It gives a Wirtinger-like presentation of any diagram group of a complete
directed $2$-complex. The translation of the proof from \cite{GuSa97} is
straightforward.

\begin{thm}
\label{th11}
Let $\kk$ be a complete directed $2$-complex with the set of negative
$2$-cells $\F^-$, and a distinguished non-empty $1$-path $w$. Then the
diagram group $\dd(\kk,w)$ admits the following presentation. The
generating set $S$ consists of all the negative edges in $\Sq(\kk,w)$
excluding edges from the left forest $T_l$. The defining relations are
all relations of the form\footnote{Here and below $x^y$ means $y^{-1}xy$.}
\be{eqf1}
(p,f_1,q\bott{f_2}r)=(p,f_1,q\topp{f_2}r)^{(\overline{p\topp{f_1}q},
f_2, r)}
\ee
if the edge $e=(\overline{p\topp{f_1}q},f_2,r)$ is not in $T_l$,
or of the form
\be{eqf2}
(p,f_1,q\bott{f_2}r)=(p,f_1,q\topp{f_2}r)
\ee
if $e\in T_l$. Here $f_1,f_2\in \F^-$, and all edges involved in
these relations are from the generating set $S$.
\end{thm}

In most cases, this presentation can be simplified.

As in \cite{GuSa97}, with every negative edge $(u,f,v)$ in
$\Sq(\kk)$ we associate a group word $[u,f,v]$ in the alphabet of
negative edges of $\Sq(\kk)$ and their inverses, defined by the
Noetherian induction on the strict order generated by $\tos$ and
the relation $\suff$, where $(u,u')\in\suff$ if and only if $u'$
is a proper suffix of $u$:

\begin{itemize}
\item If $u\ne\bar u$, then $[u,f,v]=[\bar u,f,v]$.
\item If $u=\bar u$ and $(u,f,v)$ is in $T_l$, then $[u,f,v]=1$.
\item If $u=\bar u$, $v=\bar v$ and $(u,f,v)$ is not in $T_l$,
then $[u,f,v]=(u,f,v)$.
\item If $u=\bar u$, $v\ne\bar v$ and $(u,f,v)$ is not in $T_l$,
then take the negative edge $(p,g,q)$ from the right forest $T_r$
that is assigned to $v$ (thus, $g\in \F^-$, $v=p\bott{g}q$, and
$q=\bar q$). By the induction hypothesis, we can assume that the
word $[u,f,p\topp{g}q]$ is already defined. Then let
$$
\begin{array}{l}
[u,f,v]= [u,f,p\topp{g}q]^{[\overline{u\topp{f}p},g,q]}
\end{array}\ .
$$
\end{itemize}

Notice that every letter (or its inverse) in any word $[u,f,v]$
has the form $(p,g,q)$, where $p$, $q$ are irreducible, $g\in
\F^-$, and $(p,g,q)$ is not in $T_l$.

Finally, let us present the translation of \cite[Theorem
9.8]{GuSa97} into the language of directed $2$-complexes (we are
correcting some misprints in the formulation of that theorem as
well). The translation of the proof of that theorem is
straightforward.

\begin{thm}
\label{th12}
Let $\kk$ be a complete directed $2$-complex and let $w$ be a non-empty
$1$-path in $\kk$. The group $\dd(\kk,w)$ is generated by the set $X$ of
all edges $(u,f,v)$ in $\Sq(\kk,w)$, where $u$, $v$ are irreducible,
$f\in F^-$, and $(u,f,v)$ is not in the left forest $T_l$, subject to the
following defining relations:
\be{f55}
[p,f_1,q\bott{f_2}r]=[p,f_1,q\topp{f_2}r]^{[\overline{p\topp{f_1}q},f_2,r]},
\ee
where
\begin{itemize}
\item $f_1,f_2\in F^-$,
\item $p$, $q$, $r$ are irreducible,
\item $p\topp{f_1}q\topp{f_2}r$ is homotopic to $w$ in $\kk$,
\item $(p,f_1,*)$ is not in $T_l$ and $(*,f_2,r)$ is not in $T_r$.
\end{itemize}
\end{thm}

\begin{rk}
\label{minimal}
{\rm
a) Notice that every relation in Theorem \ref{th12} is a conjugacy relation
of the form $x^y=x^z$ for some generator $x$ and words $y$, $z$. Therefore,
the set $X$ in Theorem \ref{th12} is a minimal generating set of the diagram
group (because it freely generates the abelianization of $\dd(\kk,w)$). In
Section \ref{homdg} we will show that the number of defining relations
given by Theorem \ref{th12} is also minimal possible.

b) One can check that in the formulation of Theorem \ref{th12}, we can
replace $T_r$ by $T_l$. The numbers of generators and relations will be
the same (see Remark \ref{minimal1} below), but the relations in general
will be more complicated.
}
\end{rk}

\begin{ex}
\label{ex33}
{\rm
As we mentioned before, the directed $2$-complex $\hh_n$ is complete for
every $n$. This complex has only two irreducible $1$-paths, $1$ and $x$.
Let $g_0$ be the negative $2$-cell of $\hh_n$ of the form $x=x^2$ and let
$g_1$, $g_2$, \dots, $g_n$ be the negative $2$-cells of $\hh_n$ of the form
$x=x$. Then the left forest consists of edges of the form $(1,g_0,*)^{\pm 1}$.
The right forest consists of edges of the form $(*,g_0,1)^{\pm 1}$.

By Theorem \ref{th12}, the diagram group $\dd(\hh_n,x)$ is
generated by the edges of the forms $(x,g_0,1)$, $(x,g_0,x)$, and
$(p,g_i,q)$, $1\le i\le n$, where $p,q\in \{1,x\}$ (the number of
generators is $4n+2$) subject to the conjugacy relations
(\ref{f55}), where there are $2n+1$ choices for the pair
$(p,f_1)$, $2n+1$ choices for the pair $(f_2,r)$, and $q$ can be
equal to $1$ or $x$ (the number of relations is $2(2n+1)^2$).

In particular, if $n=\infty$, the diagram group is not finitely
generated. For any integer $n$, the group $\dd(\hh_n,x)$ is finitely
presented. In case $n=0$, it has two generators $x_0=(x,g_0,1)$,
$x_1=(x,g_0,x)$ and two defining relations
$x_1^{x_0^2}=x_1^{x_0x_1}$, $x_1^{x_0^3}=x_1^{x_0^2x_1}$. This is
one of the classical presentations of R.\,Thompson's group $F$
\cite{CFP}. In case $n=1$ we get a presentation of the group
${\cal G}_1$ with $6$ generators and $18$ defining relations.
}
\end{ex}

Now we are going to use Theorems \ref{th11} and \ref{th12} to give an
example of a universal diagram group with a very simple presentation.

\begin{thm}
\label{y2y3}
Let $\kk=\la a,y\mid ay=a,y^3=y^2\ra$. Then $\kk$ is a universal
directed $2$-complex. The group $H=\dd(\kk,a)$ is universal. It
can be given by the following Thompson-like group presentation
\be{th-like}
\la x_0,x_1,x_2,\ldots\mid x_j^{x_i}=x_{j+1},\ 0\le
i<j-1\ra.
\ee
The group $H$ also has the following finite presentation with three
generators and six defining relations:
\be{th-like-fin}
\begin{array}{l}
\la x_0,x_1,x_2\mid x_2^{x_0^i}=x_2^{x_0^{i-1}x_1} (i=2,3,4),
x_2^{x_0^j}=x_2^{x_0^{j-1}x_2} (j=3,4,5)\ra.
\end{array}
\ee
\end{thm}

\proof It is easy to see that $\kk$ contains $\fft$, whence $\kk$
is universal, and $\dd(\fft,y^2)$ is embedded into $\dd(\kk,y^2)$.
Hence $\dd(\kk,y^2)$ contains copies of all countable diagram
groups. By Corollary \ref{cy2}, $\dd(\kk,a)\times\dd(\kk,y^2)$ is
embedded into $\dd(\kk,ay^2)$. Since $ay^2$ and $a$ are homotopic
in $\kk$, by Corollary \ref{cy1}, we can conclude that
$\dd(\kk,y^2)$ is embedded into $\dd(\kk,a)$. Hence $\dd(\kk,a)$
also contains copies of all countable diagram groups.

It is easy to check that the directed $2$-complex $\kk$ is complete.

Theorem \ref{th11} implies that $H$ can be generated by the edges
of the form $(u,a=ay,v)$ or $(u,y^2=y^3,v)$, where $u$, $v$ are
$1$-paths in $\kk$, $u$ is irreducible, and $uav$ (resp., $uy^2v$)
is homotopic to $a$ in $\kk$.

If $uav$ is homotopic to $a$, then clearly $u$ is empty, which
implies that $(u,a= ay, v)$ is in the left forest. Hence $H$ is
generated by the edges of the form $(u,y^2=y^3,v)$ only. If
$(u,y^2=y^3,v)$ is one of our generators, then $uy^2v$ must be
homotopic to $a$. Then $u=ay^k$, $v=y^l$ for some $k,l\ge0$. Since
$u$ must be an irreducible $1$-path, we have $k=0$. So this
generator has the form $(a,y^2=y^3,y^l)$. We denote it by $x_l$.
By Theorem \ref{th11} the defining relations of $H$ are the
following
$$
(a,y^2=y^3,vy^3w)=(a,y^2=y^3,vy^2w)^{(\overline{ay^2v},y^2=y^3,w)}.
$$
Note that $\overline{ay^2v}=a$. Let $i=|w|$, $j=|vy^2w|=i+2+|v|$.
Then our defining relation has the form $x_{j+1}=x_j^{x_i}$,
where $j\ge i+2$. This leads to (\ref{th-like}).

To describe a finite presentation of $H$, we use Theorem
\ref{th12}. Our set of generators now consists of the elements
$x_j=(a,y^2=y^3,y^j)$, where $j=0,1,2$ since the word $y^j$ is
irreducible.

The defining relations have the form
$$
[a,y^2=y^3,py^3q]=[a,y^2=y^3,py^2q]^{[a,y^2= y^3,q]}
$$
since $\overline{ay^2p}=a$, where the words $p$, $q$ are
irreducible. Also we have a restriction that $(a,y^2= y^3,q)$ is
not in the right forest. Hence $q$ is non-empty. Therefore
$p=1,y,y^2$ and $q=y,y^2$. Let $z_j=[a,y^2=y^3,y^j]$. According to
Definition \ref{dfm}, $z_j$ can be expressed as follows in terms
of the generators: $z_0=(a,y^2=y^3,1)=x_0$, $z_1=x_1$, $z_2=x_2$,
$z_3=[a,y^2= y^3,y^3]=[a,y^2=y^3,y^2]^{[a,y^2= y^3,1]} =x_2^{x_0}$
and analogously $z_4=z_3^{z_0}=x_2^{x_0^2}$, $z_5=x_2^{x_0^3}$,
and so on. Thus we have $6$ defining relations in terms of the
$z_j$'s: $z_4=z_3^{z_1}$, $z_5=z_4^{z_1}$, $z_6=z_5^{z_1}$,
$z_5=z_4^{z_2}$, $z_6=z_5^{z_2}$, $z_7=z_6^{z_2}$. If we now
rewrite them in terms of the generators $x_0$, $x_1$, $x_2$, we
obtain (\ref{th-like-fin}).
\endproof

\section{Rooted $2$-trees}
\label{2trees}

We shall need a special class of directed $2$-complexes called
{\em rooted $2$-trees\/}, which are $2$-di\-men\-si\-onal analogs of rooted
trees. Let $w$ be a simple arc subdivided into subarcs (edges). One can
regard $w$ as a directed $2$-complex with no $2$-cells. Let us denote it
by $\kk_0$. Consider any ascending family of directed $2$-complexes
\be{kkchain}
\kk_0\subseteq\kk_1\subseteq\kk_2\subseteq\cdots\subseteq\kk_n
\subseteq\cdots
\ee
obtained by the following inductive procedure. For every $n\ge0$, in order
to construct $\kk_{n+1}$, we add new $2$-cells to $\kk_n$. For each positive
$2$-cell $f$ from $\kk_{n+1}\setminus\kk_n$, the top $1$-path $\topp{f}$
must belong to $\kk_n$ and the bottom $1$-path $\bott{f}$ must be a simple
arc that meets $\kk_n$ at the endpoints only. We also assume that the arcs
of the form $\bott{f}$, for all positive $2$-cells from
$\kk_{n+1}\setminus\kk_n$, are disjoint except for their endpoints.

Then the union $\kk=\bigcup_{n\ge0}\kk_n$ is a directed $2$-complex called
a {\em rooted $2$-tree\/} (with the root $1$-path $w$).

We say that $1$-paths $p$, $q$ of a directed $2$-complex have the same
endpoints if $\iota(p)=\iota(q)$, $\tau(p)=\tau(q)$.

Let us consider the following three conditions for a directed $2$-complex
$\kk$ and a distinguished $1$-path $w$ in it:

\begin{description}
\item[T1] \ For any vertex $o$ of $\kk$, there exist a $1$-path
from $\iota(w)$ to $\tau(w)$ containing $o$.
\item[T2] \ Every two $1$-paths in $\kk$ with the same endpoints are
homotopic in $\kk$.
\item[T3] \ The diagram group $\dd(\kk,w)$ is trivial.
\end{description}

\begin{lm}
\label{2tree-prop}
Any rooted $2$-tree $\kk$ with the root $w$ satisfies conditions
{\rm T1 -- T3}.
\end{lm}

\proof
Let $\kk$ be constructed using the ascending chain (\ref{kkchain}).
It is easy to see that $\kk$ satisfies the conditions T1 -- T3 provided
each $\kk_n$ ($n\ge0$) does.

We proceed by induction on $n$. All three conditions are obvious for
$\kk_0$. So for $n\ge0$, assume that $\kk_n$ satisfies conditions
T1 -- T3.

Let us check that $\kk_{n+1}$ satisfies T1 -- T3. If $o$ is a vertex that
belongs to $\kk_{n+1}\setminus\kk_n$ then it belongs to the bottom path
$\bott{f}$ of some positive $2$-cell $f$. The endpoints of $\bott{f}$
belong to $\kk_n$ so there are $1$-paths $q'$, $q''$ in $\kk_n$ from
$\iota(w)$ to $\iota(f)$ and from $\tau(f)$ to $\tau(w)$, respectively.
Hence the $1$-path $q'\bott{f}q''$ contains $o$. Therefore, T1 holds for
$\kk_{n+1}$.

To check T2, let us consider $1$-paths $p$, $q$ in $\kk_{n+1}$, where
$\iota(p)=\iota(q)$, $\tau(p)=\tau(q)$. Let $\rho_n(s)$ denote the number
of edges in $s$ that belong to $\kk_{n+1}\setminus\kk_n$. We proceed by
induction on $\rho_n(p)+\rho_n(q)$. If this number is zero, then both
$1$-paths belong to $\kk_n$ so they are homotopic. So we can assume that
$\rho_n(q)\ne0$. This means that $q$ contains an edge $e$ from
$\kk_{n+1}\setminus\kk_n$. The edge $e$ is contained in the bottom $1$-path
$v=\bott{f}$ of a positive $2$-cell $f\not\in\kk_n$ such that $u=\topp{f}$
belongs to $\kk_n$. By definition, no edges of $v$ belong to $\kk_n$. Notice
that by definition, the graph $\kk_{n+1}$ is obtained from $\kk_n$ by adding
a union of simple arcs of the form $\bott{g}$, $g\in\kk_{n+1}\setminus\kk_n$
which are pairwise disjoint except for the endpoints. Therefore, every
$1$-path in $\kk_{n+1}$ that connects two vertices in $\kk_n$ and contains
$e$, must contain the whole $1$-path $v=\bott{f}$. Hence $q$ contains $v$.

Thus we can represent $q$ in the form $q'vq''$. Therefore, $q$ is homotopic
to the $1$-path $r=q'uq''$. Since (by definition) $u$ is contained in
$\kk_n$, we have $\rho_n(r)<\rho_n(q)$. By the inductive assumption, the
$1$-paths $p$ and $r$ are homotopic in $\kk_{n+1}$. So $p$ and $q$ are also
homotopic, that is, T2 holds for $\kk_{n+1}$.

It remains to check T3. Suppose that the group $\dd(\kk_{n+1},w)$ is not
trivial. Then there exists a reduced nontrivial $(w,w)$-diagram $\Delta$
over $\kk_{n+1}$. If all labels of the cells in $\Delta$ belong to $\kk_n$,
then $\Delta$ is a diagram over $\kk_n$ and so it represents a nontrivial
element of $\dd(\kk_n,w)$, a contradiction. So let $\pi$ be a cell in
$\Delta$ that belongs to $\kk_{n+1}$ but not to $\kk_n$. Without loss of
generality, the label $f$ of $\pi$ is a positive $2$-cell $f\in\kk_{n+1}$.
Since $\bott{f}$ is a simple arc in $\kk_{n+1}$, the path in $\Delta$
labelled by $v$ must be the common boundary of $\pi$ and a cell of the
form $\bott{f}=\topp{f}$. The only possible label of this cell is $f^{-1}$.
Thus $\Delta$ has a dipole, a contradiction.
\endproof

Now we are going to prove the converse to Lemma \ref{2tree-prop}. This
gives a characterization of rooted $2$-trees.

\begin{thm}
\label{2tree-crit}
Let $\kk$ be a directed $2$-complex and let $w$ be a $1$-path in it. Then
$\kk$ satisfies the conditions {\rm T1 -- T3} if and only if $\kk$ is a
rooted $2$-tree with the root $w$.
\end{thm}

\proof
The ``if" part is given by Lemma \ref{2tree-prop}. To prove the ``only if"
part, let us define a sequence (\ref{kkchain}) of subcomplexes in $\kk$.
Notice that by T2, any $1$-path in $\kk$ is a simple arc because an empty
$1$-path cannot be homotopic to a nonempty $1$-path. Hence we can set
$\kk_0$ to be equal to the simple arc $w$. If $\kk_n$ is already defined
for some $n\ge0$, let $\kk_{n+1}$ be the subcomplex formed by $\kk_n$ and
all the $2$-cells $f\not\in\kk_n$ such that $\topp{f}$ is contained in
$\kk_n$. Conditions T1 and T2 guarantee that $\kk$ is equal to the union of
the $\kk_n$, $n\ge0$. By definition, it is enough to show that each of the
directed $2$-complexes $\kk_n$ ($n\ge0$) is a rooted $2$-tree with $w$ as
the root.

So we prove that $\kk_n$ is a rooted $2$-tree (with $w$ as the root) by
induction on $n\ge0$. This is obvious if $n=0$. So we assume that $\kk_n$
is a rooted $2$-tree and prove the same for $\kk_{n+1}$.

In fact a stronger claim is true: for every subcomplex $\mm$ of $\kk$,
which is a rooted $2$-tree with root $w$, and any positive $2$-cell
$f\in\kk\setminus\mm$ with $\topp{f}$ in $\mm$, no internal points of
$\bott{f}$ can belong to $\mm$. (Clearly, if this claim is true, then
the directed $2$-complex $\mm'$ obtained from $\mm$ by adding $f$, is again
a rooted $2$-tree with root $w$. Since $\kk_{n+1}$ can be obtained from
$\kk_n$ by adding $2$-cells one by one, the claim implies that $\kk_{n+1}$
is a rooted $2$-tree.)

By contradiction, let $f$ be a $2$-cell in $\kk\setminus\mm$, where
$\topp{f}$ is contained in $\mm$. Suppose that there exists an internal
vertex $o$ of $\bott{f}$ that belongs to $\mm$. Since (by Lemma
\ref{2tree-prop}) $\mm$ satisfies T1, we can find a $1$-path $p$ in $\mm$
from $\iota(w)$ to $\tau(w)$ that is subdivided by $o$ into a product of
two factors, $p=p'p''$. Since the endpoints $\iota(f)$, $\tau(f)$ are also
in $\mm$, there exist $1$-paths $q'$ (respectively, $q''$) in $\mm$ from
$\iota(w)$ to $\iota(f)$ (from $\tau(f)$ to $\iota(w)$). The $1$-path
$v=\bott{f}$ is subdivided by $o$ into a product of the form $v=v'v''$. Now
we have two $1$-paths $p'$ and $q'v'$ that go from $\iota(w)$ to $o$. By T2,
they are homotopic in $\kk$. Let $\Delta'$ be a $(q'v',p)$-diagram over
$\kk$. Analogously, there exists a $(v''q'',p'')$-diagram $\Delta''$ over
$\kk$. Now we can form a $(w,w)$-diagram $\Delta$ over $\kk$ as follows. The
$1$-paths $p'p''$ and $q'uq''$ belong to $\mm$ and have the same endpoints
as $w$. Since $\mm$ satisfies T2, there exists a $(w,q'uq'')$-diagram
$\Delta_1$ and a $(p'p'',w)$-diagram $\Delta_2$. Both of them are diagrams
over $\mm$. Now let
$$
\Delta=
\Delta_1\circ(\ve(q')+\pi+\ve(q''))\circ(\Delta'+\Delta'')\circ\Delta_2,
$$
where $\pi$ is a cell of $\Delta$ labelled by $f$.

Suppose that $\Delta$ has a cell labelled by $f^{-1}$. This cell cannot
belong to $\Delta_1$ or $\Delta_2$ because they are diagrams over $\mm$.
Then this cell is in $\Delta'$ or in $\Delta''$. These cases are symmetric
so let $\pi'$ be a cell of $\Delta'$ labelled by $f^{-1}$. There is a
(directed) path in $\Delta$ from $\tau(\pi')$ to $\tau(\Delta')$. Let us
extend it by the path labelled by $v''$ from $\tau(\Delta')$ to $\tau(\pi)$.
If we project the resulting path into $\kk$, then we get a loop at $\tau(f)$.
This contradicts T2. Thus $\Delta$ contains no cells labelled by $f^{-1}$.
Therefore, after reducing all the dipoles in $\Delta$, we get a nontrivial
reduced $(w,w)$-diagram over $\kk$ which contradicts T3.
\endproof

Given a rooted $2$-tree $\kk$, the ascending sequence $\kk_n$
($n\ge0$) of rooted $2$-trees  defined in the proof of Theorem
\ref{2tree-crit}, will be called the {\em natural filtration\/} of
$\kk$. Notice that for each $2$-cell $f$ of $\kk$, there exists a
unique $n\ge1$ such that $f$ or its inverse has the form $u=v$,
where $u$ belongs to $\kk_n$ and $v$ is a simple arc in
$\kk_{n+1}\setminus\kk_n$. This allows us to choose the {\em
natural\/} orientation on the set of $2$-cells of $\kk$: we call
$f$ positive if $\topp{f}$ is contained in $\kk_n$ but $\bott{f}$
is not contained in $\kk_n$ (for some $n\ge0$).

\section{Universal $2$-covers of directed $2$-complexes}
\label{ucd2c}

Let $\kk$ be a directed $2$-complex. Consider any directed $2$-complex
$\mm$ and a morphism $\phi\colon\mm\to\kk$ that sends edges to edges
and $2$-cells to $2$-cells. For every edge or $2$-cell, its image under
$\phi$ will be called its {\em label\/}. In this situation, we shall
call $\mm$ a {\em directed $2$-complex over\/} $\kk$. The morphism $\phi$
will be sometimes called a {\em labelling map\/}. For every $1$-path $p$
in $\mm$, the set of all atomic $2$-paths in $\mm$ with top $p$ will
be called the {\em $2$-star\/} of $p$. Clearly, for any $1$-path $p$ in
$\mm$, we have an induced map from the $2$-star of $p$ into the $2$-star
of $\phi(p)$. This induced map will be called the {\em local map\/} of
$\mm$ at $p$.

Let $p$ be any $1$-path in $\kk$. A {\em universal $2$-cover\/} of
$\kk$ with base $p$ is a directed $2$-complex $\mm$ over $\kk$ with
a labelling map $\phi$ which satisfy the following properties:

\begin{description}
\item[U1]\ $\mm$ is a rooted $2$-tree with root $\tp$, where
$\phi(\tp)=p$;
\item[U2]\ for any $1$-path $\tq$ in $\mm$ from $\iota(\tp)$ to
$\tau(\tp)$, the local map of $\mm$ at $\tq$ is bijective.
\end{description}

The following theorem shows the existence and uniqueness of universal
$2$-covers. We say that two directed $2$-complexes over $\kk$ are
isomorphic whenever there exists an isomorphism between these complexes
that preserves all labels.

\begin{thm}
\label{uni-uni}
For every directed $2$-complex $\kk$ and every $1$-path $p$ in $\kk$,
there exists a universal $2$-cover of $\kk$ with base $p$. Every
two universal $2$-covers of $\kk$ with base $p$ are isomorphic.
\end{thm}

\proof
Let us construct a directed $2$-complex $\mm$ over $\kk$ defined as a
rooted $2$-tree with natural filtration
$\mm_0\subseteq\mm_1\subseteq\mm_2\subseteq\cdots$, and the labelling
map $\phi$. Both $\mm_n$ and $\phi$ are defined by induction on $n$.

By definition, $\mm_0$ is a simple arc $\tp$ subdivided into
subarcs (edges) and labelled by $p$. The labelling of $\tp$ gives
the restriction of $\phi$ onto $\mm_0$.

Suppose that a rooted $2$-tree $\mm_n$ with root $\tp$ and the
restriction $\phi$ onto $\mm_n$ are already defined for some
$n\ge0$.

Let us consider all pairs of the form $(\tr,f)$, where $\tr$ is a
$1$-path in $\mm_n$ and $f$ is a $2$-cell of $\kk$ such that
$\topp{f}=\phi(\tr)$. Suppose that $\mm_n$ does not have a
$2$-cell labelled by $f$ whose top path is $\tr$. Then we add a
new $2$-cell $\tf$ to $\mm_n$ as follows. First we add a simple
arc $\ts$ with the same endpoints as $\tr$, subdivide it into
subarcs and label it by $\bott{f}$. Then we add a $2$-cell $\tf$
of the form $\tr=\ts$, labelled by $f$. (Note that $\tf$ depends
on the pair $(\tr,f)$.) Applying this operation to all possible
pairs $(\tr,f)$ as above, we obtain $\mm_{n+1}$ and the
restriction of the labelling function $\phi$ on $\mm_{n+1}$. We
assume that all the new arcs we add are disjoint from each other
and from $\mm_n$, except for their endpoints. Clearly then,
$\mm_{n+1}$ is a rooted $2$-tree with root $\tp$.

By definition, $\mm=\bigcup_{n\ge0}\mm_n$ is also a rooted $2$-tree with
root $\tp$ so it satisfies U1. We claim that $\mm$ satisfies the
following property

\begin{description}
\item[U2$'$]\ for any $1$-path $\tr$ in $\mm$ and for any $2$-cell
$f$ of $\kk$ such that $\topp{f}=\phi(\tr)$, there is exactly one
$2$-cell $\tf$ of $\mm$ labelled by $f$ with top path $\tr$.
\end{description}

Indeed, let $n\ge0$ be the smallest integer such that $\tr$ is
contained in $\mm_n$. If there is no $2$-cell $\tf$ labelled by
$f$ with top path $\tr$ in $\mm_n$, then such a $2$-cell will
appear in $\mm_{n+1}$ by definition.

Let us prove by induction that for any $n\ge0$, there is at most
one cell $\tf$ in $\mm_n$ labelled by $f$ with the given top
$1$-path $\tr$. If $n=0$, then $\mm_0$ has no $2$-cells at all.
Suppose that our assumption is true for $\mm_n$. By the definition
of $\mm_{n+1}$, if $\mm_n$ has a $2$-cell labelled by $f$ with top
$\tr$, then no $2$-cell with the same properties appears in
$\mm_{n+1}$. Otherwise, only one such a $2$-cell may appear. Hence
our claim is true for $\mm_{n+1}$.

Now let us take any $1$-path $\tq$ in $\mm$ with the same
endpoints as $\tp$. We prove that the local map of $\mm$ at $\tq$
is surjective. Suppose that $\delta$ belongs to the $2$-star of
$q=\phi(\tq)$ in $\kk$. Then $\delta=(u,f,v)$ for some $1$-paths
$u$, $v$ in $\kk$ and some $2$-cell $f$. We also have
$q=\topp{\delta}=u\topp{f}v$. The $1$-path $\tq$ can be decomposed
as $\tq=\tu\tr\tv$, where $\phi(\tu)=u$, $\phi(\tr)=\topp{f}$,
$\phi(\tv)=v$. Since property U2$'$ holds for $\mm$, the complex
$\mm$ has a $2$-cell $\tf$ labelled by $f$ with top $\tr$. So the
atomic $2$-path $\td=(\tu,\tf,\tv)$ belongs to the $2$-star of
$\tq$ and maps onto $\delta$ under $\phi$.

Now we need to prove that the local map of $\mm$ at $\tq$ is
injective. Assume the contrary. Then there are two different
atomic $2$-paths in $\mm$ with the same top $\tq$ and the same
image $\delta=(u,f,v)$ in $\kk$. In this case the $2$-cells of
these atomic $2$-paths are also different. But they have the same
top $\tr$ and the same label $f$, which contradicts U2$'$.

We proved that the local map of $\mm$ at any $1$-path from
$\iota(\tp)$ to $\tau(\tp)$ is bijective, that is, property U2
holds. So $\mm$ is a universal $2$-cover of $\kk$ with base $p$.

Now let us prove the uniqueness. Suppose that $\mm'$ is a
universal $2$-cover of $\kk$ with base $p$. We will show that
$\mm'$ and $\mm$ are isomorphic as directed $2$-complexes over
$\kk$.

First of all, let us show that $\mm'$ also satisfies U2$'$. If $\tr$
is a $1$-path in $\mm'$ such that $\topp{f}=\phi(\tr)$ for some
$2$-cell $f$ in $\kk$, then we can include $\tr$ into a $1$-path
$\tq=\tu\tr\tv$ in $\mm$ from its $\iota(\tp)$ to $\tau(\tp)$. The
$1$-path $q=\phi(\tq)$ has the form $urv$, where
$r=\phi(\tr)=\topp{f}$. So we have an atomic $2$-path
$\delta=(u,f,v)$ in $\kk$ with $\topp{\delta}=q$. Since $\mm'$
satisfies U2, the local map of $\mm'$ at $\tq$ is surjective and
so $\delta$ has a preimage $\td$ with top $\tq$. The $2$-cell of
$\td$ is labelled by $f$ and its top is $\tr$. If there were two
$2$-cells in $\mm'$ with the same top $\tr$ and the same label
$f$, then we would have a contradiction with the injectivity of
that local map.

Let us take the natural filtration of $\mm'$, that is, the
ascending union $\mm'_0\subseteq\mm'_1\subseteq\cdots$ defined in
Section \ref{2trees}. Recall that $\mm'_0$ is the root arc of
$\mm'$ and for any $n\ge0$, $\mm'_{n+1}$ is the subcomplex in
$\mm'$ formed by $\mm'_n$ and all the $2$-cells $f\notin\mm'_n$
such that $\topp{f}$ is contained in $\mm'_n$. All the
subcomplexes $\mm'_n$ are rooted $2$-trees with the same root. It
suffices to show that $\mm_n$ and $\mm'_n$ are isomorphic as
directed $2$-complexes over $\kk$ for any $n\ge0$. We prove this
fact by induction.

If $n=0$, then both $\mm_0$ and $\mm'_0$ are simple arcs labelled
by $p$. Assume that $\mm_n$ and $\mm'_n$ are isomorphic for some
$n\ge0$ (as complexes over $\kk$). For simplicity, we can identify
$\mm_n$ and $\mm'_n$. We are going to establish a natural
bijection between the $2$-cells of $\mm_{n+1}\setminus\mm_n$ and
the $2$-cells of $\mm'_{n+1}\setminus\mm'_n$. Each pair of the
$2$-cells in the bijection, will have the same label and the same
top contained in $\mm_n=\mm'_n$. This will clearly induce the
desired isomorphism between $\mm_{n+1}$ and $\mm'_{n+1}$.

Let $\tf$ be a $2$-cell of $\mm_{n+1}\setminus\mm_n$, where
$\tr=\topp{\tf}$ is contained in $\mm_n$ and $\phi(\tf)=f$. By
U2$'$, there is a $2$-cell $\tf'$ in $\mm'$ labelled by $f$ with the
top $\tr$. This cell cannot belong to $\mm'_n$ since otherwise it
also belongs to $\mm_n$ and so $\mm$ will have two different
$2$-cells with the same label and the same top. Thus it belongs to
$\mm_{n+1}$ by definition of the natural filtration. In fact the
map $\tf\mapsto\tf'$ is the desired bijection. Indeed if $\tf'$ is
a $2$-cell of $\mm'_{n+1}\setminus\mm'_n$, then $\mm$ must contain
a $2$-cell $\tf$ labelled by $f$ with top $\tr$. By the same
arguments, $\tf$ does not belong to $\mm_n=\mm'_n$. So it belongs
to $\mm_{n+1}$ by definition.
\endproof

\begin{rk}
\label{u2pr}
{\rm
The proof of Theorem \ref{uni-uni} also shows that a directed $2$-complex
$\mm$ over $\kk$  with a labelling map $\phi$ is a universal $2$-cover of
$\kk$ if and only if it satisfies U1 and U2$'$ (thus U2 can be replaced by
U2$'$). The condition U2$'$ is sometimes easier to verify than U2.
}
\end{rk}

The universal $2$-cover of $\kk$ with base $p$ will be denoted by
$\wtk_p$. There are two useful orientations on $\wtk_p$. First, we
have the natural orientation on $\wtk_p$ because it is a rooted
$2$-tree. The set of its positive $2$-cells will be denoted by
$\wtF$. Second, given an orientation on the set of $2$-cells of
$\kk$, we have an {\em induced orientation\/} on the set of
$2$-cells of $\wtk_p$. Namely, a $2$-cell of $\wtk_p$ is positive
in that orientation whenever its image under $\phi$ is a positive
$2$-cell of $\kk$.

In order to illustrate the process of constructing $\wtk_p$, let us take
$\kk$ to be the Dunce hat $\hh_0=\la x\mid x^2=x\ra$ and let $p=x$. Then
the natural filtration $\kk_0\subseteq\kk_1\subseteq\cdots$ of $\wtk_p$
is constructed as follows. The complex $\kk_0$ is just an edge $e$ labelled
by $x$. The complex $\kk_1$ consists of three vertices, three edges,
and one positive $2$-cell (in the natural orientation) of the form
$e=e_1e_2$ labelled by the cell $x=x^2$ (all edges of $\wtk_p$ are
labelled by $x$). Now we get a new path $\tq=e_1e_2$ connecting
$\iota(e)$ and $\tau(e)$. Its image $q$ in $\hh_0$ is the $1$-path $x^2$.
There are three atomic $2$-paths in $\hh_0$ with top $1$-path $q$:
$(1,x^2=x,1)$, $(1,x=x^2,x)$ and $(x,x=x^2,1)$. The first of them
has a preimage in $\kk_1$, the other two do not have preimages.
Thus we need to add two new positive $2$-cells to $\kk_1$: a cell
of the form $e_1=e_{11}e_{12}$ and a cell of the form $e_2=e_{21}e_{22}$.
The resulting complex is $\kk_2$, it has five vertices, seven edges,
and three positive $2$-cells. The complex $\kk_3$ has seven new positive
$2$-cells: four cells of the form $e_{ij}=e_{ij1}e_{ij2}$, where
$i,j\in\{1,2\}$, and three cells $e_1e_{21}=f_1$, $e_{12}e_2=f_2$,
$e_{12}e_{21}=f$, where $f_1$, $f_2$, $f$ are new edges. Altogether the
directed $2$-complex $\kk_3$ has $9$ vertices, $18$ edges, and $10$
positive $2$-cells, etc.
\vspace{0.5ex}

In the remaining part of this section, we fix a directed $2$-complex $\kk$,
a $1$-path $p$ in $\kk$, the root $\tp$ of the universal $2$-cover $\wtk_p$
of $\kk$, the natural filtration $\kk_0\subseteq\kk_1\subseteq\cdots$ of
$\wtk_p$, and the labelling map $\phi$ from $\wtk_p$ to $\kk$.

The following lemma is the analog of the lifting lemma in the
theory of ``ordinary" covering spaces.

\begin{lm}
\label{lift}
For every $2$-path $\delta$ in $\kk$ with $\topp{\delta}=p$, there exists
a unique $2$-path $\tilde\delta$ in the universal $2$-cover $\wtk_p$ such
that $\topp{\tilde\delta}=\tp$ and $\phi(\tilde\delta)=\delta$.

For any $(p,q)$-diagram $\Delta$ over $\kk$, where $q$ is a $1$-path in
$\kk$, there exists a unique $1$-path $\tq$ in $\wtk_p$ and a unique
$(\tp,\tq)$-diagram $\wtD$ over $\wtk_p$ such that the morphism $\phi$
maps $\wtD$ onto $\Delta$.
\end{lm}

\proof
The first part of the lemma immediately follows from condition U2.

Let us consider any $(p,q)$-diagram $\Delta$ over $\kk$. It corresponds
to some $2$-path $\delta$ in $\kk$ with the top $1$-path $p$. We can lift
$\delta$ to a $2$-path $\tilde\delta$ in $\wtk_p$. The diagram $\wtD$ of
this $2$-path is a preimage of $\Delta$ under $\phi$. Clearly, any
$(\tp,\tq)$-diagram $\wtD'$ over $\wtk_p$ such that $\phi(\wtD')=\Delta$
corresponds to a lift $\td'$ of $\delta$.

All $2$-paths assigned to the same diagram $\Delta$ are isotopic
(see Section 3). Notice that the local maps preserve the property
of a pair of atomic $2$-paths to be independent. Therefore, by Lemma \ref{apis}
isotopic $2$-paths in $\kk$ lift to isotopic $2$-paths of
$\wtk_p$. Hence the diagram $\wtD$ with the desired properties is
unique.
\endproof

We shall say that the diagram $\wtD$ from the statement of Lemma \ref{lift}
is a {\em lift\/} of $\Delta$.

The following theorem shows a connection between the universal cover
of the Squier complex $\Sq(\kk)$ and the universal $2$-cover of $\kk$.

We are going to use Farley's \cite{Far00} description of the
universal cover $\Sqt(\kk)$ of $\Sq(\kk)$. Farley proved, in the case of
semigroup presentations, that $\Sqt(\kk)$ can be described in just
the same way as $\Sq(\kk)$: the vertices of $\Sqt(\kk)$ are
arbitrary diagrams over $\kk$, and for arbitrary $n\ge1$ the
$n$-cubes are pairs where the first component is an arbitrary
diagram $\Delta$ over $\kk$, and the second component is an
arbitrary thin diagram $\Theta$ over $\kk$ with $n$ cells such
that $\Delta\circ\Theta$ exists and has no dipoles. The face maps
are defined in the natural way. The restriction of the covering
map $\Sqt(\kk)\to\Sq(\kk)$ to the vertices is $\bott{\cdot}$. Any
realization of the cubical complex $\Sqt(\kk)$ will be also
denoted by $\Sqt(\kk)$. Farley's proof carries without any change
to the case of arbitrary directed $2$-complexes.

\begin{thm}
\label{ucsq}
The connected component $\Sq(\wtk_p,\tp)$ of the Squier complex of
the universal $2$-cover $\wtk_p$ is homeomorphic to the universal
cover of the connected component $\Sq(\kk,p)$ of the Squier complex
of $\kk$.
\end{thm}

\proof
We are using the fact that the space of positive paths in $\wtk_p$ is
a realization of the Squier complex $\Sq(\wtk_p)$ (Lemma \ref{OmSq}).
Consider the connected component of $\Omega_+(\wtk_p)$ that contains $\tp$.
By definition, any positive path $\tq$ in $\wtk_p$ from $\iota$ to $\tau$
can be uniquely decomposed into the product
\be{unidec}
\tq=u_0f^{(1)}_{t_1}u_1\cdots f^{(m)}_{t_m}u_m,
\ee
where $m\ge0$, $u_0,\ldots,u_m$ are $1$-paths, $f^{(1)},\ldots,f^{(m)}$
are $2$-cells in $\wtF$ and the parameters $t_1,\ldots,t_m$ belong to
the open interval $(0,1)$. It is clear that these and only these paths
form the connected component of $\tp$ in $\Omega_+(\wtk_p)$. Let
$\tr=u_0\topp{f^{(1)}}u_1\cdots\topp{f^{(m)}}u_m$. We assign to
$\tq$ an ordered pair $(\wtD,\widetilde\Psi)$, where $\wtD$ is the
reduced $(\tp,\tr)$-diagram over $\wtk_p$ (this diagram exists by
the property T2 of the rooted $2$-tree $\wtk_p$ and is unique by property T3).
By $\widetilde\Psi$ we denote the thin diagram
$\ve(u_0)+f^{(1)}+\ve(u_1)+\cdots+f^{(m)}+\ve(u_m)$. If we apply the
morphism $\phi$ to $\wtD$ and $\widetilde\Psi$, we get a pair
$(\Delta,\Psi)$, where $\Delta$ is a $(p,r)$-diagram over $\kk$ and $\Psi$
is a thin diagram over $\kk$. It is obvious that the concatenation
$\Delta\circ\Psi$ exists, moreover, it follows from the construction of
$\wtk_p$ that it has no dipoles. So the pair $(\Delta,\Psi)$ defines an
$m$-dimensional cube in $\Sqt(\kk,p)$ according to \cite{Far00}. Take
the point with coordinates $(t_1,\ldots,t_m)$ in this $m$-cube and assign
it to the element $\tq$. It is not hard to show that this defines a
homeomorphism between the connected component of the path $\tp$ in
$\Omega_+(\wtk_p)$ and the space $\Sqt(\kk,p)$ from Farley \cite{Far00}.
\endproof

The following important theorem by Farley \cite{Far00} also can be
translated into the language of directed $2$-complexes without
difficulty.

\begin{thm}
\label{Farley}
{\rm(\cite{Far00})} The universal cover $\Sqt(\kk,p)$ of $\Sq(\kk,p)$
is contractible. Hence $\Sq(\kk,p)$ is a $K(G,1)$ CW complex for the
diagram group $G=\dd(\kk,p)$, for every directed $2$-complex $\kk$
and every $1$-path $p$ in $\kk$.
\end{thm}

\begin{rk}
\label{unicov}
{\rm
The proof of contractibility of $\Sqt(\kk,p)$ given in \cite{Far00} is
based on some classical topological facts such as Whitehead's theorem.
We give a direct geometric proof using Theorem \ref{ucsq} here.

We can represent the space $\Omega_+(\wtk_p)$ as the ascending
union $\Omega_0\subseteq\Omega_1\subseteq\cdots$\,, where
$\Omega_n$ ($n\ge0$) is formed by those positive paths from the
connected component of $\tp$ which belong to $\kk_n$. Since
$\Omega_0$ consists of the single point $\tp$, it suffices to show
that $\Omega_n$ is a deformation retract of $\Omega_{n+1}$ for all
$n\ge0$. To show this, we take any element in $\Omega_{n+1}$. This
is a positive path in $\wtk_p$ from $\iota$ to $\tau$. It can be
uniquely represented in the form $h=v_0g^{(1)}_{s_1}v_1\cdots
g^{(k)}_{s_k}v_k$, where $k\ge0$, $v_0,\ldots,v_k$ are positive
paths in $\kk_n$, the $2$-cells $g^{(1)},\ldots,g^{(k)}$ belong to
$\kk_{n+1}\setminus\kk_n$ and $s_1,\ldots,s_k\in(0,1]$. Here we
also assume that all the $2$-cells $g^{(i)}$ ($1\le i\le k$) are
in $\wtF$, that is, they have their top paths in $\kk_n$ and the
bottom paths in $\kk_{n+1}\setminus\kk_n$. For any $s\in[0,1]$, we
define the positive path
\be{pphs}
h_s=v_0g^{(1)}_{ss_1}v_1\cdots g^{(k)}_{ss_k}v_k.
\ee
It is obvious that $h_0$ belongs to $\kk_n$ and $h_1=h$. It is also clear
that $h_s=h$ for all $s\in[0,1]$ provided $h$ belongs to $\kk_n$. The map
$(h,s)\mapsto h_s$ from $\Omega_{n+1}\times[0,1]$ to $\Omega_{n+1}$ is
clearly continuous. Thus $\Omega_n$ is a deformation retract of
$\Omega_{n+1}$.

The fact that $\wtk_p$ is a covering space for $\Sq(\kk,p)$ can be
proved easily as well. The covering map can be constructed as follows.
We extend $\phi$ to the set of all positive paths in $\wtk_p$ by sending
the paths of the form $\tf_t$, where $\tf$ is a $2$-cell of $\wtk_p$ and
$t\in[0,1]$, to the positive path $f_t$ in $\kk$ ($f=\phi(\tf)$ is the
label of $\tf$). Now $\phi$ maps the connected component of $\tp$ in
$\Omega_+(\wtk_p)$ onto $\Sq(\kk,p)$. We leave it to the reader as an
exercise to check that $\phi$ defines the desired covering map. (In fact,
this is essentially contained in Lemma \ref{lift} and its proof.)
}
\end{rk}

\begin{rk}
\label{FarCat}
{\rm
Recall that Farley proved that every connected component $\Sqt(\kk)$ is
a CAT(0) cubical complex provided, for example, $\kk$ is finite. Thus
is this case the diagram groups of $\kk$ act discretely, cellularly,
by isometries on CAT(0) cubical complexes.
}
\end{rk}

\begin{rk}
\label{Stall}
{\rm
There exists an alternative method of constructing the universal
$2$-cover $\wtk_p$. Similarly to the Stallings foldings in graphs
\cite{Sta83}, let us define an {\em elementary $2$-folding\/} in
a complex $\mm$ over $\kk$ as follows. If two $2$-cells $\pi_1$ and $\pi_2$ in $\mm$
have a common top path and the same label, then we identify their
bottom paths and remove one of the cells. The result of this operation
is again a directed $2$-complex over $\kk$.

Consider all possible diagrams over $\kk$ with top $1$-path labelled by
$p$. Identify the top paths of all these diagrams. Clearly, we get a
directed $2$-complex $\nn$ over $\kk$. Now do all (possibly infinitely
many) elementary $2$-foldings in $\nn$. It can be shown that the resulting
directed $2$-complex over $\kk$ is $\wtk_p$. We shall prove this in our
next paper where foldings are used in a more general situation.
}
\end{rk}

\section{Homology}
\label{homdg}

Theorem \ref{Farley} shows that the components $\Sq(\kk,w)$ are
$K(G,1)$ spaces for diagram groups $G=\dd(\kk,w)$. In fact in most
cases $\Sq(\kk)$ is too large. Here we will use the technique of
collapsing schemes \cite{BrGe,Br92,Coh} to find a ``smaller" CW
complex, which is homotopy equivalent to $\Sq(\kk)$ (at least in
the case when $\kk$ is complete).

We recall the concept of collapsing scheme from
\cite{Br92,Coh,BrGe}. Let $X$ be a semi-cubical complex. We say
that we have a collapsing scheme for $X$ if the following is true:

\begin{itemize}
\item there exists a subdivision of the set of all cubes of $X$
into three disjoint subsets: {\em essential\/}, {\em collapsible\/},
and {\em redundant\/} cubes;
\item there exists a strict partial order $\succ$ on the set of
all redundant $n$-cubes ($n\ge0)$ that satisfies the descending
chain condition (that is, any sequence $c_1\succ c_2\succ\cdots$
terminates).
\item there exists a bijection $c\mapsto\hat c$ between the set of all
redundant $n$-cubes and the set of all collapsible $(n+1)$-cubes (for
every $n\ge0$);
\item any redundant $n$-cube $c$ occurs exactly once among the
$n$-faces of $\hat c$ and all the other redundant $n$-faces $c'$ of
$\hat c$ precede $c$ in the order $\succ$ (that is, $c\succ c'$); the
redundant $n$-cube $c$ is called the {\em free face\/} of the
collapsible $(n+1)$-cube $\hat c$.
\end{itemize}

The next lemma is proved in almost the same way as
\cite[Proposition 1]{Br92}.

\begin{lm}
\label{Br-coll}
Given a collapsible scheme for a semi-cubical complex $X$, one can
construct a CW complex $Y$ which is homotopy equivalent to $X$ in such a
way that the $n$-cells of Y are in one-to-one cor\-res\-pon\-den\-ce with
the essential $n$-cubes of $X$.
\end{lm}

As in \cite{Br92}, for each $n\ge0$ one has to do an infinite
number of elementary steps, one for each collapsible $n$-cube. The
free face of a collapsible cube is identified (homeomorphically)
with the union of the other faces and the collapsible cube
disappears. The space $Y$ is a end result of the process.

Now let $\kk$ be a complete directed $2$-complex. Let $T_l$ be a
left forest in $\kk$. Recall that for any $1$-path $p$, the
irreducible form of $p$ is denoted by $\bar p$.

Let $c=\ve(u_0)+f_1+\cdots+f_n+\ve(u_n)$ be an $n$-cube of $X$
(a thin diagram with $n$ cells). For any $0\le i\le n$, we say
that the term $u_i$ in $c$ is {\em special\/} provided it is not
an irreducible $1$-path, that is, $\bar u_i\ne u_i$. For $1\le i\le n$,
we say that the term $f_i$ in $c$ is {\em special\/} provided
$(u_{i-1},f_i,*)$ is in $T_l$.

The $n$-cube $c$ is called {\em essential\/} if it has no special
terms. If $c$ is not essential, then we find its leftmost special
term. If this is one of the $u_i$'s, then we call $c$ {\em redundant\/}.
Otherwise we call $c$ {\em collapsible\/} (in this case the special term
is one of the $f_i$'s).

To describe the strict partial order $\succ$, we need to introduce one
technical concept. Let $\Psi=\ve(u_0)+g_1+\cdots+g_k+\ve(u_k)$
($k\ge0$) be a thin diagram over $\kk$. Suppose that for some
$i$, $j$, where $0\le i\le j\le k$, we have $u_j=u'\bott{g}u''$,
for some $1$-paths $u'$, $u''$ and $g\in\F^-$. Then it is possible to
{\em move the cell $g_i$ to the right\/} replacing $\Psi$ by a new thin
diagram $\Psi'=\ve(u_0)+\cdots+\ve(u_{i-1})+\ve(\bott{g_i})+\cdots+\ve(u')+g+
\ve(u'')+\cdots$ (we replace the term $g_i$ by $\ve(\bott{g_i})$ and the
term $\ve(u_j)$ by $\ve(u')+g+\ve(u'')$, thus we remove the cell $g_i$ and
insert a cell $g$). It is easy to see that this process of moving cells to
the right always terminates. Indeed, $\Psi'$ has also $k$ cells and it can
be represented in the form $\Psi'=\ve(u_0')+g_1'+\cdots+g_k'+\ve(u_k')$. If
we compare the $(k+1)$-vectors ${\bf b}=(|u_k|,\ldots,|u_0|)$ and
${\bf b}'=(|u_k'|,\ldots,|u_0'|)$, then it follows from our description
that ${\bf b}'$ strictly precedes ${\bf b}$ in the lexicographical
order. Since these vectors have non-negative coordinates, the process of
moving cells to the right must terminate.

Now we can define $\succ$. Let $c$, $c'$ be redundant $n$-cubes of
$X$. If $\bott{c}\tos\bott{c'}$ then we set $c\succ c'$. Otherwise,
if $\bott{c}=\bott{c'}$, then we set $c\succ c'$ whenever $c'$ can
be obtained from $c$ by a (non-zero) number of moving cells to the
right. The fact that $\kk$ is Noetherian and the remark from the
previous paragraph imply that $\succ$ is a strict partial order
satisfying the descending chain condition.

Let $c=\ve(u_0)+f_1+\cdots+f_n+\ve(u_n)$ be a redundant $n$-cube.
This means that $u_i$ is not irreducible for some $0\le i\le n$
whereas all the terms to the left of $u_i$ are not special. Let us
find the edge $(p,f,q)$ in $T_l$ assigned to $u_i$. Thus
$u_i=p\bott{f}q$. Note that $p$ is irreducible by definition. Thus
the $(n+1)$-cube $\ve(u_0)+f_1+\cdots+f_i+\ve(p)+f+\ve(q)+\cdots$
is collapsible. We denote it by $\hat c$ and check that $c$ is the
free face of $\hat c$. We consider all $n$-faces $\topp{\hat c}_j$
and $\bott{\hat c}_j$ of $\hat c$, $1\le j\le n+1$. Let $j=i+1$.
Then $\bott{\hat c}_j=c$ and $c'=\topp{\hat c}_j$ is obtained from
$c$ by replacing $\bott{f}$ by $\topp{f}$. Since $f$ participates
in an edge from $T_l$, one has $\bott{f}\tos\topp{f}$. Hence
$c\succ c'$ whenever $c'$ is redundant.

Now suppose that $j>i+1$. Then $\bott{\hat c}_j$ and $\topp{\hat c}_j$ are
collapsible $n$-cubes. Thus we may skip this case because we compare $c$
with redundant cubes only.

Finally, take $j\le i$. We can only consider the case when
$c'=\bott{\hat c}_j$ since $\bott{f_j}\tos\topp{f_j}$ or
$\bott{f_j}=\topp{f_j}$. Suppose that $c'$ is a redundant $n$-cube.
To compare $c$ and $c'$, notice that their bottom paths are the same.
The diagram $c'$ is obtained by moving one cell of $c$ to the right
(the cell $f_j$ has been deleted and the cell $f$ has been added). Thus
$c\succ c'$.

It remains to check that we have a bijection between redundant
$n$-cubes and collapsible $(n+1)$-cubes. We already assigned a
collapsible $(n+1)$-cube $\hat c$ to each redundant $n$-cube $c$.
If we start with a collapsible $(n+1)$-cube, then we can find its
leftmost special term. This is some $f\in\F^{-}$. The cube then has
the form $\cdots+\ve(p)+f+\ve(q)+\cdots$, where $(p,f,q)$ is in $T_l$.
Replacing $f$ by $\bott{f}$ gives us a redundant cube $c$. It follows
directly from definitions that the cube $\hat c$ assigned to $c$ is
exactly the collapsible $(n+1)$-cube we started with. This completes our
proof that we have defined a collapsible scheme for $X$.

Summarizing and taking into account Lemma \ref{Br-coll}, we get
the following.

\begin{thm}
\label{KG1}
Suppose that $\kk$ is a complete directed $2$-complex and let
$G=\dd(\kk,w)$, where $w$ is a non-empty $1$-path in $\kk$. Then there
exists a $K(G,1)$ CW complex $Y_w$ whose $n$-dimensional cells are in
one-to-one correspondence with thin diagrams of the form
$c=\ve(u_0)+f_1+\cdots+f_n+\ve(u_n)$, where $n\ge0$, the $1$-paths $u_i$
are irreducible for all $0\le i\le n$, the edges $(u_{i-1},f_i,*)$ are not
in $T_l$ for all $1\le i\le n$, and $\topp{c}$ is homotopic to $w$ in $\kk$.
\end{thm}

Note that each thin diagram $c$ described in the statement of
Theorem \ref{KG1} is an essential cube in the Squier complex
$\Sq(\kk)$.

Recall that a directed $2$-complex is called $2$-path connected if
all non-empty $1$-paths in it are homotopic. Let us call a directed
$2$-complex {\em almost $2$-path connected\/} if the number of
classes of homotopic $1$-paths is finite (that is, $\Sq(\kk)$ has
finitely many connected components). A complex of the form $\kk_\pp$,
where $\pp$ is a semigroup presentation, is almost $2$-path connected
if and only if the semigroup given by $\pp$ is finite. Thus the
following statement generalizes a result from \cite{Far00} and
strengthens \cite[Theorem 10.7]{GuSa97}.

\begin{thm}
\label{KGfinite}
Let $\kk$ be a finite almost $2$-path connected $2$-complex. Then all
diagram groups of $\kk$ are of type $\fff_\infty$.
\end{thm}

\proof
First suppose that $\kk$ is complete. Let $C$ denote the number of homotopy
classes of $1$-paths in $\kk$ (the empty $1$-paths are included) and let $N$
be the number of positive $2$-cells of $\kk$. It is clear that the number of
$n$-cells in the $K(G,1)$ space $Y_w$ does not exceed $C(CN)^n$.
In particular, it is finite so $G$ has type $\fff_\infty$.

Now suppose that $\kk$ is not necessarily complete. By Lemma
\ref{retract}, $\kk$ is contained  in a finite complete almost $2$-path
connected directed $2$-complex $\kk'$, and the diagram groups
of $\kk$ are retracts of the diagram groups of $\kk'$. It remains
to recall that a retract of an $\fff_\infty$ group is of type
$\fff_\infty$.
\endproof

By \cite[Theorem 10.3]{GuSa97}, if $\pp$ is a finite complete
rewrite system such that all diagram groups over it are finitely
generated, then all of them are finitely presented. Now we can
deduce a much stronger result. In fact we can even eliminate the
assumption that the presentation is finite.

\begin{thm}
\label{fp1finfty}
Let $\kk$ be a complete directed $2$-complex. Suppose that all diagram
groups of $\kk$ are finitely generated. Then all of them are of type
$\fff_\infty$.
\end{thm}

\proof
By Theorem \ref{KG1} it is enough to prove that for any $n\ge1$, each
connected component of $\Sq(\kk)$ has only finitely many essential
$n$-cubes. We proceed by induction on $n$. Let $n=1$. By definition,
the set of essential cubes of dimension $1$ is in one-to-one correspondence
with the generating set of the corresponding diagram group described in
Theorem \ref{th12}. Since this set is minimal by Remark \ref{minimal},
it is finite.

Now let $n>1$. For any essential $1$-cube $\ve(p)+f+\ve(q)$ of a
connected component $\Sq(\kk,w)$ of $\Sq(\kk)$, let us consider
the set of all essential $(n-1)$-cubes in $\Sq(\kk,q)$. By the
inductive assumption, it is finite. It remains to note that any
essential $n$-cube $c=\ve(u_0)+f_1+\cdots+f_n+\ve(u_n)$ of
$\Sq(\kk,w)$ is determined uniquely by an essential $1$-cube
$\ve(u_0)+f_1+\ve(q)$ of $\Sq(\kk,w)$ and an essential
$(n-1)$-cube $\ve(u_1)+f_2+\cdots+f_n+\ve(u_n)$ from $\Sq(\kk,q)$.
\endproof

\begin{rk}
\label{strng}
{\rm
Note that one can extract a stronger fact from the proof of Theorem
\ref{fp1finfty}. Suppose that $\kk$ is a complete directed $2$-complex.
If some diagram group $\dd(\kk,w)$ is not of type $\fff_\infty$, then
there are two $1$-paths $w'w_1$ and $w_2w''$ homotopic to $w$ such that
the diagram groups $\dd(\kk,w')$ and $\dd(\kk,w'')$ are not finitely
generated.
}
\end{rk}

Now we are going to prove that for any complete directed
$2$-complex $\kk$, the complex $Y_w$ described in the statement
of Theorem \ref{KG1} is in fact ``minimal". Namely, for every
$n\ge 0$,  we shall compute the integer $n$th homology group of
every diagram group $\dd(\kk,w)$ of $\kk$ and show that it is a
free Abelian group whose rank coincides with the number of
$n$-cells in $Y_w$ (that is, the number of essential $n$-cubes in
$\Sq(\kk,w)$.

Let $G=\dd(\kk,w)$. Since the homology groups of a group $G$
coincide with the homology groups of any $K(G,1)$ CW complex, let
us consider the complex $X=\Sq(\kk,w)$ (it is a $K(G,1)$ by
Theorem \ref{Farley}). As usual, let $T_l$ be a left forest in $X$.

Denote by $P_n$ the free Abelian group with the set of $n$-cubes
of $X$ as a free basis. The boundary maps $\partial_n\colon P_n\to
P_{n-1}$ ($n\ge1$) are given by the formulas of Serre \cite[p.\,440]{Ser51}
(see also \cite{BrGe}):
\be{bou}
\partial_n(c)=\sum\limits_{i=1}^n (-1)^i(\topp{c}_i-\bott{c}_i),
\ee
where $c$ is an $n$-cube. Since the maps (\ref{bou}) form a
chain complex \cite{Ser51}, the $n$th integer homology group
$H_n(G;\zz)$ coincides with the $n$th homology group of that
chain complex (that is, $\Ker\partial_n/\Img\partial_{n+1}$).

As in \cite{Br92,Coh}, we define an endomorphism $\Ev$ of the chain
complex $P=(P_n,\partial_n)$. This endomorphism maps
every $P_n$ into the subgroup $Q_n$ of $P_n$ freely generated by
the essential $n$-cubes. Let $c$ be an $n$-cube (a generator of
$P_n$). If $c$ is collapsible then we set $\Ev(c)=0$. If $c$ is
essential then we set $\Ev(c)=c$. Finally suppose that $c$ is
redundant. In that case we proceed by the Noetherian induction on
$\succ$.

Since $c$ is redundant, there exists a collapsible $(n+1)$-cube
$\hat c$ such that $c$ is the free face of $\hat c$. Then
$\partial_{n+1}(\hat c)=\pm c+\Sigma$ for some linear combination
$\Sigma$ of cubes that are either essential or collapsible or
redundant but smaller than $c$ with respect to $\succ$. Thus we
can assume that $\Ev(c')$ has been defined already for all $c'$
occurring in $\Sigma$. So we can set $\Ev(c)=\mp\Ev(\Sigma)$.

It is shown in \cite{Coh} (see also \cite[p.\,150]{Br92}), that $\Ev$
indeed is an endomorphism of the chain complex (that is, it
commutes with the boundary maps), and that the chain complex $Q$
formed by the groups $Q_n$ and boundary maps
$\delta_n=\Ev\partial_n$ is chain-equivalent to the initial chain
complex. Thus the homology groups of $Q$ coincide with the homology
groups of $P$.

We are going to prove that $\delta_n$ is a zero map. For that we
need the following statement.

\begin{lm}
\label{rewrpr}
For any $n$-cube $c=\ve(u_0)+f_1+\cdots+f_n+\ve(u_n)$ of $\Sq(\kk,w)$,
we let $\bar c=\ve(\bar u_0)+f_1+\cdots+f_n+\ve(\bar u_n)$. Then
$\Ev(c)=\bar c$ if $\bar c$ is essential and $\Ev(c)=0$ if $\bar c$ is
collapsible $($note that $\bar c$ cannot be redundant by definition$)$.
\end{lm}

\proof
Note that if $c$ is essential then $c=\bar c$ and $\Ev(c)=c$ as required.

Suppose that $c$ is collapsible. Then $\Ev(c)=0$ by definition.
Thus we only need to check that $\bar c$ is also collapsible. Since
$c$ is collapsible, the edge $(u_{i-1},f_i,*)$ is in the left forest
$T_l$ for some $i\le n$ and all the $u_j$'s are irreducible
for $0\le j<i$. Then $\bar c$ possesses similar properties whence
$\bar c$ is collapsible.

Now suppose that $c$ is redundant. Take the smallest number $i\le n$
such that $u_i$ is not irreducible, and consider the edge
$(p,f,q)$ from $T_l$ assigned to $u_i$. This $i$ will be called
the {\em index\/} of $c$. Consider also the collapsible $(n+1)$-cube
$\hat c=\ve(u_0)+\cdots+f_i+\ve(p)+f+\ve(q)+\cdots$ whose free
cell is $c$. Since $\kk$ is Noetherian, we can assume without loss
of generality that the statement of the lemma does not hold for
$c$ but holds for all $n$-cubes $c'$ such that $\bott{c}\tos\bott{c'}$.
We can also assume that among all such counterexamples, $c$ has the smallest
index $i$.

Notice that for every $j$, $1\le j\le n+1$,
\be{nnn}
\overline{\topp{\hat c}_j}=\overline{\bott{\hat c}_j}.
\ee

Also notice that $c=\bott{\hat c}_{i+1}$ and $\bott{c}\tos\bott{c'}$,
where $c'=\topp{\hat c}_{i+1}$. By (\ref{nnn}), $\bar c=\overline{c'}$.
By the definition of $\Ev$, we have $\Ev(c)=\Ev(c')+\Ev(\Sigma)$, where
$\Sigma$ is the sum of $\topp{\hat c}_j-\bott{\hat c}_j$, $j\ne i+1$ by
(\ref{bou}).

It remains to check that
\be{equal}
\Ev(\topp{\hat c}_j-\bott{\hat c}_j)=0
\ee
for every $j\ne i+1$. If $j>i+1$, then both cells $\topp{\hat c}_j$ and
$\bott{\hat c}_j$ are collapsible, so
$\Ev(\topp{\hat c}_j)=\Ev(\bott{\hat c}_j)=0$ and (\ref{equal}) holds.

Let $1\le j\le i$. The thin diagrams $e'=\topp{\hat c}_j$ and
$e=\bott{\hat c}_j$ are obtained from $\hat c$ by replacing $f_j$
by its top and bottom path, respectively. Suppose that $e'\ne e$
(otherwise there is nothing to prove). Thus $\bott{f_j}\tos\topp{f_j}$ and
the statement of the lemma holds for $e'$ because
$\bott{c}=\bott{e}\tos\bott{e'}$. It is easy to see from definition that
the cube $e$ is redundant. The index of $e$ is $j-1<i$. Hence the statement
of the lemma holds for $e$ as well. By (\ref{nnn}), $\overline{e'}=\bar e$.
Therefore, $\Ev(e'-e)=0$ and (\ref{equal}) holds.
\endproof

Now it is easy to show that all boundary maps $\delta_n$, $n\ge1$,
in the chain complex $Q$ are zero. Indeed, by Lemma \ref{rewrpr},
the value $\Ev(c)$ depends only on $\bar c$. By (\ref{nnn}) and
(\ref{bou}), for every $n$-cube $c$,
$$
\delta_n(c)=\Ev\partial_n(c)=\Ev\left(\sum\limits_{i=1}^n
(-1)^i(\topp{c}_i-\bott{c}_i)\right)=\sum\limits_{i=1}^n
(-1)^i(\Ev(\topp{c}_i)-\Ev(\bott{c}_i))=0.
$$

Thus $\Ker\partial_n=Q_n$ and $\Img\partial_{n+1}=0$. Hence for
$H_n(G;\zz)\cong Q_n$ is free Abelian, $n\ge1$.

If $n=0$, then $H_0(G;\zz)\cong\zz$ and we have only one essential
cell of dimension $0$ --- this is the vertex corresponding to the
irreducible $1$-path of $\Sq(\kk,w)$. Thus we proved

\begin{thm}
\label{genhom}
Let $\kk$ be a complete directed $2$-complex and let $w$ be a non-empty
$1$-path in $\kk$. The $n$th integer homology $H_n(G;\zz)$ $(n\ge0)$ of
the diagram group $G=\dd(\kk,w)$ is free Abelian. Its free basis consists
of all essential $n$-cubes from $\Sq(\kk,w)$.
\end{thm}

Theorem \ref{genhom} implies that the CW complex $Y_w$ from
Theorem \ref{KG1} gives a minimal presentation of $G=\dd(\kk,w)$
in terms both the number of generators and the number of relations.
In fact, it gives a minimal set of generators of homology groups in all
dimensions.

\begin{rk}
\label{minimal1}
{\rm
It is not difficult to prove that the presentation given by $Y_w$
is precisely the presentation from Theorem \ref{th12}, where $T_r$
is replaced by $T_l$ (see Remark \ref{minimal}, part  b). We
already know (Remark \ref{minimal}, part a) that the presentation
from Theorem \ref{th12} involves the minimal possible number of
generators. Let us show that it contains the minimal number of
relations as well. For any $1$-path $p$, let us denote by $\mu(p)$
the minimal number of generators for the diagram group $G=\dd(\kk,w)$.
This is exactly the number of edges of the form $(u,f,v)\notin T_l$
that belong to $\Sq(\kk,w)$, where $u$, $v$ are irreducible,
$f\in F^{-}$ (these are the essential $1$-cubes of $\Sq(\kk,w)$).
If we replace here $T_l$ by $T_r$, then we again have a minimal
generating set of $G$ because of a symmetry. It is easy to give a
formula to compute the number of the defining relations of $G$ given
by Theorem \ref{KG1}. Let $s_1$, \dots, $s_m$ be the third components
of the essential $1$-cubes of $G$ ($m=\mu(w)$). Each of the defining
relations corresponds to an essential $1$-cube in $\Sq(\kk,s_i)$ for
some $i$. Thus the sum $\mu(s_1)+\cdots+\mu(s_m)$ is exactly the number
of the defining relations of $G$ given by Theorem \ref{KG1}. Clearly,
it will be the same if we use $T_r$ instead of $T_l$ in the definition
of essential cubes. Thus the number of defining relations given by
Theorem \ref{th12} and Theorem \ref{KG1} are the same.
}
\end{rk}

Now consider the homology groups of arbitrary diagram groups. Let $H$ be a
diagram group of an arbitrary directed $2$-complex. We know from Lemma
\ref{emb1}, part 4 that $H$ is a retract of a diagram group $G$ of a
complete directed $2$-complex. Notice that the retraction can be described
in the language of group homomorphisms: $H$ is a retract of $G$ if and only
if there are two homomorphisms $\phi\colon G\to H$ and $\psi\colon H\to G$
such that $\phi\psi={\bf id}_H$ ($\psi$ acts first). Since $H_n(-,\zz)$
is a covariant functor \cite{Br-Cohom}, this implies that $H_n(H;\zz)$
is a retract of $H_n(G;\zz)$. In particular, $H_n(H;\zz)$ is also free
Abelian and its rank does not exceed the rank of $H_n(G;\zz)$. So we proved

\begin{thm}
\label{homfree}
For any $n\ge0$ and for any diagram group $G$, the $n$th integer homology
group $H_n(G;\zz)$ is free Abelian.
\end{thm}

This theorem answers a question by S. Pride.
\vspace{1ex}

If $G$ is a group of type $\fff_\infty$, then one can consider its
Poincar\'{e} series
$$
P_G(t)=\sum\limits_{n=0}^\infty r_nt^n,
$$
where $r_n$ denotes the rank of the $n$th integer homology group
of $G$. Note that $r_0=1$.

\begin{ex}
\label{abc}
{\rm
a) Let $\kk=\la x\mid x^r=x\ra$ ($r\ge2$). It is proved in \cite{GuSa97}
that $\dd(\kk,x)\cong F_r$, where $F_r$ is a generalization of the
R.\,Thompson group $F=F_2$ defined in \cite{Br92}. That complex is complete
so we can use Theorems \ref{KG1} and \ref{genhom}. The essential $n$-cells
of this complex ($n\ge1$) have the form
$c=\ve(x^{k_0})+(x=x^r)+\ve(x^{k_1})+\cdots+(x=x^r)+\ve(x^{k_n})$,
where $1\le k_i<r$ for $0\le i<n$, $0\le k_n<r$. These cells may
belong to different components of $\Sq(\kk)$. Clearly, $c$ belongs
to the component of $x$ if and only if the sum
$n+k_0+k_1+\cdots+k_n$ equals $1$ modulo $r-1$. This condition
determines uniquely the number $k_0$ by the other numbers. So
there are exactly $r\cdot(r-1)^{n-1}$ ways to choose $c$ in
$\Sq(\kk,x)$. (For instance, in the case $r=2$ we have the
R.\,Thompson's group $F$, which has $\zz^2$ as its $n$th integer
homology group for $n\ge1$, which was proved in \cite{BrGe}.)

The Poincar\'{e} series for $F_r$ has the form
$$
P(t)=1+rt+r(r-1)t^2+\cdots+r(r-1)^{n-1}t^n+\cdots=\frac{1+t}{1-(r-1)t}.
$$

b) Now let $\fft=\la y\mid y^3=y^2\ra$ be the complex on Figure 9
(by Theorem \ref{th24}, the diagram group $\dd(\fft,y^2)$ is
universal). The complex $\fft$ is complete as well. The essential
$n$-cubes in $\Sq(\fft,x^2)$ have the form
$c=\ve(y^{k_0})+(y^2=y^3)+y^{k_1}+\cdots+(y^2=y^3)+\ve(y^{k_n})$,
where $k_i=1,2$ for $0\le i<n$, $k_n=1,2,3$. All of them for
$n\ge1$ belong to $\Sq(\kk,y^2)$. Hence the rank of the $n$th
homology group of $\dd(\fft,x^2)$ is $3\cdot2^n$ ($n\ge1$) and so
the Poincar\'{e} series is
$$
P(t)=1+6t+12t^2+\cdots+3\cdot2^nt^n+\cdots=\frac{1+4t}{1-2t}.
$$
Notice that the Poincar\'e series of $\dd(\fft,x^2)$ coincides
with the Poincar\'e series of the free product $F_3*F_3$ but these
diagram groups are not isomorphic (which can be proved by using
Kurosh's theorem).

c) One more universal diagram group is given by the complete
directed $2$-complex $\hh_1=\la x\mid x^2=x,x=x\ra$ (Theorem \ref{th565}).
Let us find the number of the essential $n$-cubes in $\Sq(\hh_1,x)$.
Let $c=\ve(u_0)+f_1+\cdots+f_n+\ve(u_n)$ be one of these cubes. Then
there are three possibilities for each of the pairs $(u_{i-1},f_i)$,
namely, $(x,x=x^2)$, $(1,x=x)$, or $(x,x=x)$ ($1\le i\le n$). There are
two possibilities for the $1$-path $u_n$, namely, $1$ and $x$. So the
rank of the $n$th homology group of ${\cal G}_1=\dd(\hh_1,x)$ equals
$2\cdot3^n$ for $n\ge1$. It is interesting to mention that both universal
groups considered in b), c) have the same minimal number of generators
(equal to $6$) and we know that they are embeddable into each other
because of their universal property. However, they are not isomorphic
because their second homology groups have ranks $12$ and $18$, respectively.
Thus the Poincar\'{e} series of ${\cal G}_1$ is
$$
P(t)=1+6t+18t^2+\cdots+2\cdot3^nt^n+\cdots=\frac{1+3t}{1-3t}.
$$
}
\end{ex}

We see that all these Poincar\'{e} series are rational. This can
be explained by the following

\begin{thm}
\label{ratio}
Let $\kk$ be a complete finite almost $2$-path connected directed
$2$-complex. Then the Poincar\'{e} series of any of its diagram groups
is rational.
\end{thm}

\proof
We refer to \cite{Lall} for the well-known properties of rational
languages. Let $A=\PP\cup(\PP\times \F^-)$, where $\PP$ consists of
all irreducible $1$-paths in $\kk$, including the empty $1$-paths,
$\F^-$ consists of all negative $2$-cells of $\kk$.

Let $0$ be a symbol not in $\PP$, and let us define a binary operation
$\cdot$ on $S=\PP\cup\{0\}$: $p\cdot q=\overline{pq}$ if $\tau(p)=\iota(q)$,
$p,q\in \PP$; all other products are equal to 0. It is easy to see that $S$
is a semigroup.

Let $\phi$ be a homomorphism from the free semigroup $A^+$ to $S$ induced
by the map that takes each $p\in \PP$ to itself and each pair
$(p,f)\in\PP\times \F^-$ to $\overline{p\topp{f}}=\overline{p\bott{f}}$.
Notice that for every $s\in \PP$, every word $w$ of the form
$(u_0,f_1)\cdots(u_{n-1},f_n)u_n$ from the rational language
$\phi^{-1}(s)\subseteq A^+$ corresponds to an $n$-cube in $\Sq(\kk,s)$. This
cube is essential if and only if $w$ does not contain letters of the form
$(p,f)$, where $p\in\PP$, $f\in\F^-$, and $(p,f,*)$ belongs to the fixed
left forest $T_l$ of $\Sq(\kk)$.

Thus let $L$ be the sublanguage of $\phi^{-1}(s)$ consisting of all
words from $\phi^{-1}(s)\cap A^+\PP$ which do not contain the letters
described in the previous paragraph. Clearly, $L$ is a rational language.
There exists a one-to-one correspondence between words of length $n+1$ in
$L$ and essential $n$-cubes in $\Sq(\kk,s)$. Since the generating function
of a rational language $L$ is rational (see for example \cite{ChSch}), the
Poincar\'{e} series of $\Sq(\kk,s)$ is a rational function as well.
\endproof

Recall \cite{Br-Cohom} that for any group $G$, its geometric
dimension, $\gd(G)$, is the smallest dimension of a $K(G,1)$. Its
cohomological dimension, $\cd(G)$, is the length of the shortest
projective resolution of the trivial $\zz G$-module $\zz$. It is
easy to see \cite{Br-Cohom} that
\be{cdgd}
\cd(G)\le\gd(G).
\ee
By the Eilenberg -- Ganea theorem \cite{LS77}, $\cd(G)=\gd(G)$ provided
$\cd(G)\ne2$ or $\gd(G)\ne3$.

Theorems \ref{KG1} and \ref{genhom} immediately imply that for
diagram groups over complete directed $2$-complexes these two
dimensions coincide.

\begin{thm}
\label{cd}
For every diagram group $G$ over a complete directed $2$-complex,
$$
\cd(G)=\gd(G).
$$
\end{thm}

\proof
Let $G=\dd(\kk,w)$. By Theorem \ref{genhom} the length of any projective
resolution of the trivial $\zz G$-module $\zz$ cannot be smaller than the
highest dimension $n$ of an essential cube of $\Sq(\kk,w)$. By Theorem
\ref{KG1}, $n\ge\gd(G)$. Therefore, $\cd(G)\ge n\ge\gd(G)$. Hence by
(\ref{cdgd}), $\gd(G)=\cd(G)$.
\endproof

The following result gives an algebraic characterization of groups
of finite cohomological dimension among diagram group of complete
directed $2$-complexes.

\begin{thm}
\label{cd1}
Let $G$ be a diagram group over a complete directed $2$-complex $\kk$,
and $n$ be a natural number. Then $\cd(G)\ge n$ if and only if $G$
contains a copy of $\zz^n$.
\end{thm}

\proof
Let $G=\dd(\kk,w)$. Let $T_l$ be a left forest in $\Sq(\kk)$. The ``if"
statement is well known \cite{Br-Cohom}. So suppose that $\cd(G)\ge n$.
Then by Theorem \ref{genhom}, $\Sq(\kk,w)$ contains an essential cube
$c=\ve(u_0)+f_1+\cdots+f_n+\ve(u_n)$.

By definition, for every $1\le i\le n$, the edge $(u_{i-1},f_i,*)$ does
not belong to $T_l$. Hence the connected component
$\Sq(\kk,u_{i-1}\bott{f_i})$ contains edges not in $T_l$. Therefore, by
Theorem \ref{th12} (or Theorem \ref{genhom}), the group
$G_i=\dd(\kk,u_{i-1}\bott{f_i})$ is non-trivial. Since all diagram groups
are torsion-free \cite{GuSa97}, $G_i$ contains a copy of $\zz$.

By Corollary \ref{cy2}, the diagram group $G_1\times\cdots\times G_n$ embeds
into $H=\dd(\kk,p)$, where $p=u_0\bott{f_1}\cdots\bott{f_n}u_n$. Therefore,
a copy of $\zz^n$ is contained in $H$. But by the choice of $c$, the
$1$-paths $p$ and $w$ are in the same connected component of $\Sq(\kk)$.
Hence by Corollary \ref{cy1}, $H$ is isomorphic to $G$, so $G$ contains a
copy of $\zz^n$, as required.
\endproof

Notice that Theorem \ref{comm} below gives a characterization of
directed $2$-complexes $\kk$ such that $\dd(\kk,w)$ contains a copy
of $\zz^n$.

Theorem \ref{cd1} immediately implies the following result.

\begin{thm}
\label{cd2}
A diagram group $G$ over a complete directed $2$-complex is free if and
only if $G$ does not contain a copy of $\zz^2$. In particular, a hyperbolic
group can be a diagram group of a complete directed $2$-complex if and only
if it is free.
\end{thm}

\proof
Indeed, by Theorem \ref{cd}, if $G$ does not contain $\zz^2$ then
$\cd(G)\le1$, and one can use the well known result of Stallings and Swan
\cite{LS77} (or, easier, one can use Remark \ref{minimal1}, and conclude
that $G$ has a presentation with no relations).
\endproof

\begin{prob}
Is it possible to drop the completeness restriction from the formulations
of Theorems {\rm \ref{ratio}, \ref{cd}, \ref{cd1}, \ref{cd2}}?
\end{prob}

\begin{rk}
{\rm
Recall that by Lemma \ref{OmSq} the space of positive paths $\Omega_+(\kk)$
is a realization of $\Sq(\kk)$. By Theorem \ref{Farley}, the homology of a
connected component of that space coincides with the homology of the corresponding
diagram group. Hence by Theorem \ref{ratio}, the Poincar\'e series of the space of
positive paths of a complete almost $2$-path connected directed $2$-complex is
rational. This resembles the well known result of Serre (see, for example, \cite{An})
that the Poincar\'e series of the loop space of a simply connected CW $2$-complex is
always rational.
}
\end{rk}

\begin{rk}
\label{ababab}
{\rm
Notice that the  completeness restriction in the statements of this paper
can be replaced by the condition ``there exists a left forest". Say, let
$\kk=\la a,b\mid ab=a,ba=b\ra$. It is not hard to check that $\kk$ is not
complete. However, one can construct a spanning forest satisfying conditions
F1 and F2 of the left forest (it is formed by all edges of the form
$(1,a=ab,*)$, $(a,b=ba,*)$, $(1,b=ba,*)$, $(b,a=ab,*)$ and the inverse
edges). Using that forest, as above, one can compute the presentation of
the corresponding diagram groups, and their homology groups. The Poincar\'e
series of the diagram groups of this complex are rational.
}
\end{rk}

\section{Rigidity}
\label{rigidity}

Recall that the flat torus theorem \cite[Theorem 7.1]{Bri} says, in
particular, that if $X$ is a metric space with CAT(0) universal cover
$\tilde X$ and the fundamental group of $X$ contains a copy of
$\mathbb{Z}^n$, then $X$ contains a $\pi_1$-embedded torus
$\mathbb{R}^n/\zz^n$.

Results of this section are of similar spirit: we prove that a
diagram groupoid of a directed $2$-complex $\kk$ contains certain
diagram group $G=\dd(\sss,p)$ if and only if there exists a
$p$-nonsingular morphism from $\sss$ into $\kk$. Since every
$p$-nonsingular morphism $\phi\colon\sss\to\kk$ induces a
$\pi_1$-injective continuous map $\Sq(\sss,p)\to\Sq(\kk,\phi(p))$,
these results (and Theorem \ref{fg}) imply that if
$\pi_1(\Sq(\kk),u)$ contains a copy of $G=\pi_1(\Sq(\sss),p)$
then there exists a $\pi_1$-injective continuous map from
$\Sq(\sss,w)$ into $\Sq(\kk)$.

In general we say that a triple (diagram group $G$, directed
$2$-complex $\kk$, $1$-path $p$ in $\kk$) is {\em rigid\/} if
$G=\dd(\kk,p)$ and for every directed complex $\kk'$ such that
$\dd(\kk')$ contains a copy of $G$ there exists a $p$-nonsingular
morphism of $\kk$ into $\kk'$.

\begin{center}
\unitlength=1mm
\special{em:linewidth 0.4pt}
\linethickness{0.4pt}
\begin{picture}(26.67,24.67)
\put(1.00,8.67){\line(1,0){25.67}}
\bezier{124}(1.00,8.67)(15.00,17.33)(26.67,8.67)
\bezier{120}(1.00,8.67)(15.67,1.00)(26.67,8.67)
\bezier{216}(1.00,8.67)(14.33,32.67)(26.67,8.67)
\put(13.67,2.33){\makebox(0,0)[cc]{$x$}}
\put(13.67,6.67){\makebox(0,0)[cc]{$x$}}
\put(13.67,11.00){\makebox(0,0)[cc]{$x$}}
\put(13.67,16.00){\makebox(0,0)[cc]{$\ldots$}}
\put(13.67,22.67){\makebox(0,0)[cc]{$x$}}
\put(1.00,5.67){\makebox(0,0)[cc]{$\iota$}}
\put(26.67,6.33){\makebox(0,0)[cc]{$\tau$}}
\end{picture}

\nopagebreak[4] Figure \theppp.
\end{center}

For example, consider the directed $2$-complex $\kk$ on Figure
\theppp\ with two vertices $\iota$ and $\tau$, one edge $x$
connecting $\iota$ with $\tau$ and positive $2$-cells of the form
$x=x$ labelled by elements of some set $A$. By Theorem \ref{th12}
(or a straightforward computation), the diagram group $G=D(\kk,x)$
is the free group of rank $|A|$. It is easy to see (exercise) that
the triple $(G,\kk,x)$ is rigid.

\addtocounter{ppp}{1}

A much more nontrivial example of a rigid triple involves the
R.\,Thompson group $F$. Example \ref{ex33} shows that $F$ is the
diagram group of the Dunce hat $\hh_0$.

The following proposition shows a remarkable property of the Dunce hat.

\begin{thm}
\label{h0tok}
Every morphism from the Dunce hat to any directed $2$-complex ${\kk}$
is nonsingular.
\end{thm}

\proof
Let us consider any morphism $\phi$ from $\hh_0$ into a directed
$2$-complex $\kk$.  We need to show that it is $p$-nonsingular
for every non-empty $1$-path $p$ in $\hh_0$. Since for every such
$p$, the diagram groups $\dd(\hh_0,p)$ and $\dd(\hh_0,x^5)$ are
conjugate, it is enough to show that $\phi$ is $x^5$-nonsingular.
Since $x^5$ and $x$ are homotopic, we have
$\dd(\hh_0,x^5)\cong\dd(\hh_0,x)\cong F$ (see Example \ref{ex33}).

Since all proper homomorphic images of the group $F$ are Abelian
\cite{CFP}, it suffices to find  an  element in the derived subgroup of
$\dd(\hh_0,x^5)$ that is not in the kernel of $\phi_{x^5}$. Let $\Pi$ be
the diagram over $\hh_0$ corresponding to the atomic $2$-path $(1,f,1)$,
where $f$ is the positive $2$-cell $x^2=x$ in $\hh_0$. Consider the diagram
$\Psi=\ve(x)+\Pi+\Pi^{-1}+\ve(x)$ from $\dd(\hh_0,x^5)$. By Theorem 11.3
from \cite{GuSa97}, $\Psi$ is in the derived subgroup of $\dd(\hh_0,x^5)$
and is nontrivial. By Theorem \ref{fg}, we can assume that the diagram
$\Delta=\phi_{x^5}(\Pi)$ is reduced. Since it is a
$(\phi(x)^2,\phi(x))$-diagram, it is nontrivial. Then
$\phi_{x^5}(\Psi) =\ve(\phi(x))+\Delta+\Delta^{-1}+\ve(\phi(x))$ is
also a reduced and nontrivial diagram, hence by Theorem \ref{fg},
$\phi_{x^5}(\Psi)\ne 1$.
\endproof

Note that the same property is true for directed $2$-complexes $\la
x\mid x^r=x\ra$ that correspond to groups $F_r$ ($r\ge2)$, the
generalizations of $F$ (see \cite{Bro87,GuSa97}). The proof is
based on the same idea (all proper homomorphic images of these
groups are also Abelian).

Here is a reformulation of the main result of \cite{GuSa02} which
shows that the triple $(F,\hh_0,x)$ is rigid.

\begin{thm}
\label{th2}
{\rm(\cite{GuSa02})} Let $\kk$ be a directed complex. Then the following
conditions are equivalent.
\begin{enumerate}
\item A diagram groupoid $\kk$ contains an isomorphic copy of the
R.\,Thompson group $F$.
\item The complex $\kk$ contains a non-empty $1$-path which is
homotopic to its square.
\item There exists a $($nonsingular$)$ morphism from the Dunce hat to
$\kk$.
\end{enumerate}
\end{thm}

Thus if the diagram groupoid of $\kk$ contains a copy of $F$ then it
contains a naturally embedded copy of $F$.

Another example of a rigid triple is given by \cite[Theorem 24]{GuSa99}.
Let $\qq$ be the directed $2$-complex with three vertices, three edges
$x$, $y$, $z$ and three positive cells of the forms $xy=x$, $y=y$, $yz=z$
on Figure \theppp\ (to obtain the complex from the diagram , we identify
all edges having the same labels).

It is proved in \cite{GuSa99} that the diagram group $\dd(\qq,xyz)$ is
isomorphic to the restricted wreath product $\zz\wr\zz$. Theorem 24 of
\cite{GuSa99} shows that the triple $(\zz\wr\zz,\qq,xyz)$ is rigid.

\begin{center}
\unitlength=1.00mm
\special{em:linewidth 0.4pt}
\linethickness{0.4pt}
\begin{picture}(90.33,35.67)
\put(3.33,17.34){\line(1,0){29.33}}
\put(61.33,17.34){\line(1,0){27.33}}
\put(3.67,17.34){\circle*{1.33}}
\put(32.33,17.34){\circle*{1.33}}
\put(61.67,17.34){\circle*{1.33}}
\put(89.67,17.34){\circle*{1.33}}
\bezier{360}(3.67,17.34)(30.33,51.67)(61.67,17.34)
\bezier{340}(32.33,17.34)(64.33,-14.00)(89.67,17.34)
\bezier{164}(32.33,17.34)(45.00,31.67)(61.67,17.34)
\bezier{164}(32.33,17.34)(45.00,4.67)(61.67,17.34)
\put(32.00,32.67){\makebox(0,0)[cc]{$x$}}
\put(18.33,14.67){\makebox(0,0)[cc]{$x$}}
\put(47.33,12.67){\makebox(0,0)[cc]{$y$}}
\put(74.67,14.67){\makebox(0,0)[cc]{$z$}}
\put(45.33,21.67){\makebox(0,0)[cc]{$y$}}
\put(63.00,4.34){\makebox(0,0)[cc]{$z$}}
\end{picture}

\nopagebreak[4] Figure \theppp.
\end{center}
\addtocounter{ppp}{1}

The free Abelian group  $\zz^n$ can participate in a rigid triple
too (for every $n\ge1$). In fact, using a description of commuting
diagrams (\cite[Theorem 17]{GuSa99}) one can obtain a much more precise
result (Theorem \ref{comm} below).

Let $\sss_n$ be the following directed $2$-complex: take a simple
path labelled by the word $w_n=x_1\cdots x_n$, where $x_i$ are
letters and let us attach $n$ positive $2$-cells $x_1=x_1$, \dots,
$x_n=x_n$ to it. Thus $\sss_n$ is a chain of $n$ spheres (Figure \theppp).

\begin{center}
\unitlength=1mm
\special{em:linewidth 0.4pt}
\linethickness{0.4pt}
\begin{picture}(115.00,27.00)
\bezier{168}(5.67,16.33)(20.00,31.00)(35.33,16.33)
\bezier{164}(5.33,16.33)(20.67,2.00)(35.33,16.33)
\bezier{168}(35.67,16.33)(50.00,31.00)(65.33,16.33)
\bezier{164}(35.33,16.33)(50.67,2.00)(65.33,16.33)
\bezier{168}(85.33,16.33)(99.67,31.00)(115.00,16.33)
\bezier{164}(85.00,16.33)(100.33,2.00)(115.00,16.33)
\put(75.00,16.67){\makebox(0,0)[cc]{$\ldots$}}
\put(20.33,25.67){\makebox(0,0)[cc]{$x_1$}}
\put(20.33,5.67){\makebox(0,0)[cc]{$x_1$}}
\put(50.33,5.67){\makebox(0,0)[cc]{$x_2$}}
\put(50.33,25.67){\makebox(0,0)[cc]{$x_2$}}
\put(100.00,25.67){\makebox(0,0)[cc]{$x_n$}}
\put(100.00,5.67){\makebox(0,0)[cc]{$x_n$}}
\end{picture}

\nopagebreak[4] Figure \theppp.
\end{center}
\addtocounter{ppp}{1}

\noindent
Then $\Sq(\sss_n,w_n)$ is the $n$-dimensional torus (this is easy
to check). Thus by Theorem \ref{fg}, the group $\dd(\sss_n,w_n)$
is isomorphic to $\zz^n$.

The next result is one of the most useful technical facts about diagram
groups. In \cite{GuSa99}, it is formulated and proved for diagrams over
semigroup presentations (see \cite[Theorem 24]{GuSa99}). The proof for
directed complexes is completely similar (in fact it can be deduced from
the result of \cite{GuSa99} by using subdivisions of complexes). It is
similar to the well known theorem that commuting matrices over an
algebraically closed field are simultaneously conjugate to their Jordan
forms.

\begin{lm}
\label{prws-cmmt}
Let $\kk$ be a directed $2$-complex and $w$ be a non-empty $1$-path in
$\kk$. Suppose that $A_1$, \dots, $A_n$ are spherical $(w,w)$-diagrams
that pairwise commute in $G$. Then there exist a $1$-path $v=v_1\ldots v_m$,
spherical $(v_j,v_j)$-diagrams $\Delta_j$ $(1\le j\le m)$ over $\kk$,
integers $d_{ij}$ $(1\le i\le n$, $1\le j\le m)$ and some $(w,v)$-diagram
$\Gamma$ over $\kk$ such that
$$
\Gamma^{-1}A_i\Gamma=\Delta_1^{d_{i1}}+\cdots+\Delta_m^{d_{im}}
$$
for all $1\le i\le n$.
\end{lm}

\begin{thm}
\label{comm}
The triple $(\zz^n,\sss_n,w_n)$ is rigid for every $n\ge1$. In addition,
let $\kk$ be a directed $2$-complex. Then every copy of $\zz^n$ in
$\dd(\kk)$ is conjugate in $\dd(\kk)$ to a subgroup of a naturally
embedded copy of $\zz^m$ for some $m\ge n$.
\end{thm}

\proof
Suppose that a diagram group $\dd(\kk,p)$ of some directed $2$-complex
contains a copy $G=\la A_1,\ldots,A_n\ra$ of $\zz^n$. We use the notation
from Lemma \ref{prws-cmmt}. Let us also assume that $m$ is chosen to be
minimal. Then all diagrams $\Delta_1$, \dots, $\Delta_m$ are nontrivial.
Indeed, if we assume the contrary, then $m>1$ because $G\ne1$. If
$\Delta_i$ is trivial for some $i$, then one has $i<m$ or $i>1$. Without
loss of generality we assume that $i<m$. But now it would be possible to
replace $\Delta_i$ by $\Delta_{i+1}$ taking into account that the power
of a sum is the sum of powers and all powers of a trivial diagram coincide.

Clearly, $m\ge n$ (otherwise the rank of the subgroup generated by
$A_1$, \dots, $A_n$ would be less than $n$).

Now the map $\phi$ from $\sss_m$ to $\kk$ that sends the positive $2$-cell
$x_i=x_i$ to the $2$-path corresponding to $\Delta_i$ ($1\le i\le m$),
defines a morphism. The image of $\dd(\sss_m,w_m)$ under $\phi_{w_n}$
is generated by the diagrams
$\ve(v_1\cdots v_{i-1})+\Delta_i+\ve(v_{i+1}\cdots v_m)$. Thus the image
of $\phi_{w_m}$ is isomorphic to $\zz^m$. Therefore, $\phi$ is
$w_m$-nonsingular. Clearly, the naturally embedded copy
$\phi(\dd(\sss_m,w_m))$ of $\zz^m$ contains $\Gamma^{-1}G\Gamma$. This
proves the second statement of the theorem.

In order to prove the rigidity statement, we just note that $\sss_n$ is a
subcomplex of $\sss_m$ so it maps into it nonsingularly. But we already
have a nonsingular morphism of $\sss_m$ into $\kk$. It suffices to compose
the morphisms.
\endproof

\begin{prob}
{\rm
It is interesting to characterize other diagram groups that can
participate in rigid triples. In particular, in view of rigidity of
the triple $(F,\hh_0,x)$ it is natural to ask if the analog is true
for $\hh_n$. By Theorem \ref{th565}, it is enough to prove that for
$n=1$. This would give a characterization of universal directed
$2$-complexes as those admitting a nonsingular morphism from $\hh_1$.
}
\end{prob}

\vspace{10ex}

\begin{minipage}[t]{3in}
\noindent Victor Guba\\
Department of Mathematics\\
Vologda State University\\
guba@uni-vologda.ac.ru
\end{minipage}
\begin{minipage}[t]{3in}
\noindent Mark V. Sapir\\
Department of Mathematics\\
Vanderbilt University\\
msapir@math.vanderbilt.edu
\end{minipage}

\end{document}